\documentclass[12pt]{amsart} 
\usepackage{amsfonts,amsmath,amsthm,amssymb,latexsym,amsthm,newlfont,enumerate,color}

\newif\ifslide

\theoremstyle{plain}

\ifdefined\spacer
\newtheorem{theorem}{Theorem}
 \else
\newtheorem{theorem}{Theorem}[section]
\fi

\newtheorem{corollary}[theorem]{Corollary}
\newtheorem{lemma}[theorem]{Lemma}

\newtheorem*{theorem*}{Theorem}
\newtheorem*{corollary*}{Corollary}

\newtheorem{proposition}[theorem]{Proposition}
\newtheorem{definition-lemma}[theorem]{Definition-Lemma}
\newtheorem{question}[theorem]{Question}
\newtheorem{red-question}[theorem]{\textcolor{red}{Question}}

\theoremstyle{definition}
\newtheorem{definition}[theorem]{Definition}
\newtheorem{remark}[theorem]{Remark}

\newtheorem{example}[theorem]{Example}

\newcommand{\dto}{\dashrightarrow}

\def\ideal#1.{I_{#1}}
\def\ring#1.{\mathcal {O}_{#1}}

\def\fring#1.{\hat{\mathcal {O}}_{#1}}
\def\proj#1.{\mathbb {P}(#1)}
\def\pr #1.{\mathbb {P}^{#1}}
\def\dpr #1.{\hat{\mathbb {P}}^{#1}}
\def\af #1.{\mathbb A^{#1}}
\def\Hz #1.{\mathbb F_{#1}}
\def\Hbz #1.{\overline{\mathbb F}_{#1}}
\def\fb#1.{\underset #1 {\times}}
\def\rest#1.{\underset {\ \ring #1.} \to \otimes}
\def\au#1.{\operatorname {Aut}\,(#1)}
\def\deg#1.{\operatorname {deg } (#1)}

\def\pic#1.{\operatorname {Pic}\,(#1)}
\def\pico#1.{\operatorname{Pic}^0(#1)}
\def\picg#1.{\operatorname {Pic}^G(#1)}
\def\ner#1.{NS (#1)}
\def\rdown#1.{\llcorner#1\lrcorner}
\def\rfdown#1.{\lfloor{#1}\rfloor}
\def\rup#1.{\ulcorner{#1}\urcorner}
\def\rcup#1.{\lceil{#1}\rceil}

\def\n1#1.{\operatorname {N_1}(#1)}  
\def\cn1#1.{\overline{\operatorname {N^1}(#1)}} 
\def\cone#1.{\operatorname {NE}(#1)}     
\def\ccone#1.{\overline{\operatorname {NE}}(#1)}
\def\none#1.{\operatorname {NF}(#1)}
\def\cnone#1.{\overline{\operatorname {NF}}(#1)}
\def\mone#1.{\operatorname {NM}(#1)} 
\def\cmone#1.{\overline{\operatorname {NM}}(#1)}

\def\coef#1.{\frac{(#1-1)}{#1}}
\def\vit#1.{D_{\langle #1 \rangle}}
\def\mm#1.{\overline {M}_{0,#1}}
\def\H1#1.{H^1(#1,{\ring #1.})}
\def\ac#1.{\overline {\mathbb F}_{#1}}

\def\adj#1.{\frac {#1-1}{#1}}
\def\spn#1.{\overline{#1}}
\def\pek#1.#2.{\Cal P^{#1}(#2)}
\def\plk#1.#2.{\Cal P^{\leq #1}(#2)}
\def\ev#1.{\operatorname{ev_{#1}}}
\def\ilist#1.{{#1}_1,{#1}_2,\dots}
\def\bminv#1.{(\nu_1,s_1;\nu_2,s_2;\dots ;\nu_{#1},s_{#1};\nu_{r+1})}
\def\zinv#1.{(\nu_1,s_1;\nu_2,s_2;\dots ;\nu_{#1},s_{#1};0)}
\def\iinv#1.{(\nu_1,s_1;\nu_2,s_2;\dots ;\nu_{#1},s_{#1};\infty)}

\def\scr #1.{\mathcal #1}

\def\llist#1.#2.{{#1}_1,{#1}_2,\dots,{#1}_{#2}}
\def\ulist#1.#2.{{#1}^1,{#1}^2,\dots,{#1}^{#2}}
\def\lomitlist#1.#2.{{#1}_1,{#1}_2,\dots,\hat {{#1}_i}, \dots, {#1}_{#2}}
\def\lomitlistz#1.#2.{{#1}_0,{#1}_1,\dots,\hat {{#1}_i}, \dots, {#1}_{#2}}
\def\loc#1.#2.{\Cal O_{#1,#2}}
\def\fderiv#1.#2.{\frac {\partial #1}{\partial #2}}
\def\deriv#1.#2.{\frac {d #1}{d #2}}

\def\map#1.#2.{#1 \longrightarrow #2}
\def\rmap#1.#2.{#1 \dasharrow #2}
\def\emb#1.#2.{#1 \hookrightarrow #2}
\def\non#1.#2.{\text {Spec }#1[\epsilon]/(\epsilon)^{#2}}
\def\Hi#1.#2.{\text {Hilb}^{#1}(#2)}
\def\sym#1.#2.{\operatorname {Sym}^{#1}(#2)}
\def\Hb#1.#2.{\text {Hilb}_{#1}(#2)}
\def\Hm#1.#2.{\Hom_{#1}(#2)}
\def\prd#1.#2.{{#1}_1\cdot {#1}_2\cdots {#1}_{#2}}
\def\Bl #1.#2.{\operatorname {Bl}_{#1}#2}
\def\pl #1.#2.{#1^{\otimes #2}}
\def\mgn#1.#2.{\overline {M}_{#1,#2}}
\def\ialist#1.#2.{{#1}_1 #2 {#1}_2, #2\dots}
\def\pair#1.#2.{\langle #1, #2\rangle}
\def\vandermonde#1.#2.{\left|
\begin{matrix}
1 & 1 & 1 & \dots & 1\\
{#1}_1 & {#1}_2 & {#1}_3 & \dots & {#1}_{#2}\\
{#1}_1^2 & {#1}_2^2 & {#1}_3^2 & \dots & {#1}_{#2}^2\\
\vdots & \vdots & \vdots & \ddots & \vdots\\
{#1}_1^{#2-1} & {#1}_2^{#2-1} & {#1}_2^{#2-1} & \dots & {#1}_{#2}^{#2-1}\\
\end{matrix}
\right|
}
\def\vandermondet#1.#2.{\left|
\begin{matrix}
1 & {#1}_1   & {#1}_1^2 & \dots & {#1}_1^{#2-1}\\
1 & {#1}_2   & {#1}_2^2 & \dots & {#1}_2^{#2-1}\\
1 & {#1}_3   & {#1}_3^2 & \dots & {#1}_3^{#2-1}\\
\vdots & \vdots & \vdots & \ddots & \vdots\\
1 & {#1}_{#2}& {#1}_{#2}^2 & \dots & {#1}_{#2}^{#2-1}\\
\end{matrix}
\right|
}
\def\gr#1.#2.{\mathbb{G}(#1,#2)}

\def\alist#1.#2.#3.{{#1}_1 #2 {#1}_2 #2\dots #2 {#1}_{#3}}
\def\zlist#1.#2.#3.{#1_0 #2 #1_1 #2\dots #2 #1_{#3}}
\def\lomitlist30#1.#2.#3.{{#1}_0,{#1}_1 #2 \dots #2\hat {{#1}_i} #2\dots #2 {#1}_{#3}}
\def\lmap#1.#2.#3.{#1 \overset{#2}{\longrightarrow} #3}
\def\mes#1.#2.#3.{#1 \longrightarrow #2 \longrightarrow #3}
\def\ses#1.#2.#3.{0\longrightarrow #1 \longrightarrow #2 \longrightarrow #3 \longrightarrow 0}
\def\les#1.#2.#3.{0\longrightarrow #1 \longrightarrow #2 \longrightarrow #3}
\def\res#1.#2.#3.{#1 \longrightarrow #2 \longrightarrow #3\longrightarrow 0}
\def\Hi#1.#2.#3.{\text {Hilb}^{#1}_{#2}(#3)}
\def\ten#1.#2.#3.{#1\underset {#2}{\otimes} #3}
\def\lomitlist30#1.#2.#3.{{#1}_0 #2 {#1}_1 #2 \dots #2 \hat {{#1}_i} #2 \dots #2 {#1}_{#3}}
\def\mderiv#1.#2.#3.{\frac {d^{#3} #1}{d #2^{#3}}}

\def\Hom{\operatorname{Hom}}

\def\Proj{\operatorname{Proj}}

\def\dim{\operatorname{dim}}

\def\deg{\operatorname{deg}}

\def\GL{\operatorname{GL}}
\def\SL{\operatorname{SL}}

\def\det{\operatorname{det}}

\def\Sing{\operatorname{Sing}}

\def\rk{\operatorname{rk}}

\def\rest{\operatorname{res}}

\def\vol{\operatorname{vol}}

\def\SL{\operatorname{SL}} 

\def\C{\mathbb C}

\def\e{\Cal E}
\def\ZZ{\mathbb Z}

\def\e1{E_1}
\def\e2{E_2}

\def\Q{\mathbb Q}
\def\Z{\mathbb Z}

\def\mapdown#1{\big\downarrow\rlap{$\vcenter{\hbox{$\scriptstyle#1$}}$}}

\def\mapse#1{
{\vcenter{\hbox{$\mathop{\smash{\raise1pt\hbox{$\diagdown$}\!\lower7pt
\hbox{$\searrow$}}\vphantom{p}}\limits_{#1}\vphantom{\mapdown{}}$}}}}

\def\VR#1.{height#1pt&\omit&&\omit&&\omit&&\omit&&\omit&\cr}

\def\VRT#1.{height#1pt&\omit&&\omit&\cr}

\usepackage[colorlinks=true,pagebackref,hyperindex]{hyperref}

\title{On the Chern numbers of a smooth threefold}

\author{Paolo Cascini}
\address{Department of Mathematics\\
Imperial College London\\
180 Queen's Gate\\
London SW7 2AZ, UK}
\email{p.cascini@imperial.ac.uk}
\author{Luca Tasin}
\address{Universit\`a Roma Tre, Dipartimento di Matematica e Fisica, Largo San Leonardo Murialdo I-00146 Roma, Italy} 
\email{ltasin@mat.uniroma3.it}

\thanks{The first author was partially supported by an EPSRC Grant. The second author was supported by the DFG Emmy Noether-Nachwuchsgruppe ``Gute Strukturen in der h\"oherdimensionalen birationalen Geometrie''. The second author was funded by the Max Planck Institute for Mathematics in Bonn during part of the realization of this project. We would like to thank G. Codogni,  S. Lu, C. Mautner, M. McQuillan, D. Panov and S. Schreieder for several very useful discussions.
We would also like to thank the referee for carefully reading our manuscript and for suggesting several improvements.
}

\begin{document}

\begin{abstract}We study the behaviour of Chern numbers of  three dimensional terminal varieties under divisorial contractions. 
\end{abstract}
\maketitle
\tableofcontents

\section{Introduction}\label{s_introduction}

The main goal of this paper is to study the Chern numbers of a smooth projective threefold, especially in relation with divisorial contractions. To this aim we will investigate the interplay between topological properties and birational properties of 3-folds.

The starting point of our research is the following question of Hirzebruch \cite{Hirzebruch54}: Which linear combinations of Chern numbers on a smooth complex projective variety are topologically invariant? 

Hirzerbruch's question has been answered by Kotschick  \cite{Kotschick08, Kotschick12}, who showed that 
a rational linear combination of Chern numbers is a homeomorphism invariant of smooth complex projective varieties if and only if it is a multiple of the Euler characteristic. 
In particular, Kotschick shows the existence of a sequence of infinitely many pairs of smooth projective threefolds $X_i, Y_i$, with $i\in \mathbb N$, such that $X_i$ and $Y_i$ are diffeomorphic and
$$c_1c_2(X_i)\neq c_1c_2(Y_i)\qquad \text{and}\qquad c_1^3(X_i)\neq c_1^3(Y_i)$$
for each $i\in \mathbb N$.

In view of this, it is natural to ask if the Chern numbers of an $n$-dimensional smooth projective variety can only assume finitely many values, after we fix the underlying manifold. 
In general, $c_n$ is a topological invariant, as it coincides with the Euler characteristic, and therefore if $n=1$ then the problem is easily settled. On the other hand, if $X$ and $Y$ are homeomorphic complex surfaces, then either $c_1^2(X)=c_1^2(Y)$ or $c_1^2(X)= 4c_2(Y)-c_1^2(Y)$, depending on whether the homeomorphism between $X$ and $Y$ is orientation preserving or not (cfr. \cite{Kotschick08}). Nevertheless, if $X$ and $Y$ are diffeomorphic surfaces, then $c_1(X)^2=c_1(Y)^2$.

In dimension three, the relevant Chern numbers are $c_1c_2$ and $c_1^3$. If $X$ is K{\"a}hler, then by  the Hirzebruch-Riemann-Roch theorem we have
$$
\vert \frac{1}{24}c_1c_2(X) \vert=| \chi(\mathcal O_X)|=| 1-h^{1,0}+h^{2,0}-h^{3,0} | \le 1+b_1 + b_2+b_3,
$$ 
where $h^{i,0}=h^i(X,\ring X.)$ and $b_1,b_2$ and $b_3$ denote the topological Betti numbers of $X$. Thus, $c_1 c_2(X)$ is bounded by a linear combination of the Betti numbers of $X$. On the other hand, LeBrun \cite{Le99} shows that the same result does not hold if we drop the assumption of being K\"ahler, answering a question raised by Okonek and Van de Ven \cite{OV95}. In  particular, he shows that if $M$ denotes the $4$-manifold underlying a K3 surface and $S^2$ is the two dimensional sphere, then there exist infinitely many complex structures $J_m$ on $M\times S^2$  such that $c_1c_2=48m$, with $m\in \mathbb N$.

More generally, in  dimension $n$,  Libgober and  Wood \cite{LW90} showed that $c_1c_{n-1}$ can be expressed in terms of Hodge numbers  and, in particular,  it is bounded by a constant that depends only on the Betti numbers of the underlying topological space. Recently, Schreieder and the second author  \cite{ST15} studied the problem in dimension at least $4$, proving that in complex dimension $n \ge 4$, the Chern numbers $c_n$, $c_1c_{n-1}$ and and $c_2^2$ ($n = 4$) are the only Chern numbers that take on only finitely many values on
the complex projective structures with the same underlying smooth 2n-manifold.

\medskip

Thus, the motivating question of this paper is  the following 

\begin{question}\cite[Problem 1]{Kotschick08}\label{q}
Does $c_1^3=-K_X^3$  take only finitely many values on the projective algebraic structures $X$ with the same underlying 6-manifold?
\end{question}

Our aim is to study this problem from a birational point of view. 

\medskip

Let $X$ be a smooth threefold. We first consider Question \ref{q} in three extreme cases which arise as building blocks in birational geometry: Fano manifolds, Calabi-Yau and canonically polarized varieties. In the first case,  it is known that $X$ belongs to a bounded family and in particular $K_X^3$ is bounded \cite{Kollar93b}. If $X$ is a Calabi-Yau, then by definition $K_X=0$ and therefore $K_X^3=0$. 
Finally, if $X$ is canonically polarized (i.e. $K_X$ is ample), then the Bogomolov-Miyaoka-Yau inequality implies that $0<K_X^3\le 8/3 c_1c_2(X)$. Thus,  the arguments above imply that $K_X^3$ is bounded by the Betti numbers of $X$.

We now consider the general case of a smooth projective threefold $X$.  Thanks to Mori's program \cite{KM98}, we can run a Minimal Model Program (MMP, in short) on $X$ and obtain  a birational map $\varphi\colon X\dashrightarrow Y$ into a threefold $Y$
such that either $X$ is not uniruled and $Y$ is minimal (i.e. the canonical divisor $K_Y$ is nef) or $X$ is uniruled and $Y$ admits a Mori fibre space structure (i.e. a  morphism $Y\to Z$ with connected fibres with relative Picard number equal to one and whose general fibre is a non-trivial Fano variety).  Thus, our strategy consists  in two steps: we first want to bound $K_Y^3$ and then  bound $K_X^3-K_Y^3$. 

One of the difficulties of the first step is due to the fact that  in general $Y$ is not smooth, but it admits some mild singularities, called terminal. On the other hand, by \cite{CZ14}, we can bound the singularities of $Y$, and in particular the index of each singularity, by a bound which depends only on the topology of $X$ (see Proposition \ref{p_2b2}). 

Recall that a variety of dimension $n$ is said to be \textit{uniruled} if there exists a variety $Y$ of dimension $n-1$ and a dominant rational map $Y\times \mathbb P^1\dashrightarrow X$. In particular, if $X$ is uniruled then it is covered by rational curves, i.e. for each $x\in X$ there exists a non-trivial morphism $f\colon \mathbb P^1\to X$ such that $x\in f(\mathbb P^1)$.

Note that if $X$ is not uniruled then $Y$ is minimal and $K_Y^3$ coincides with the volume of $X$ (cf. definition \ref{d_volume}), which is a birational invariant of the variety $X$. 

Our first result, based on Bogomolov-Miyaoka-Yau inequality for terminal threefolds, is the following: 

\begin{theorem}\label{t_vol}
Let $X$ be a smooth complex projective threefold which is not uniruled. 

Then 
$$\vol(X,K_X)\le 6b_2(X)+36b_3(X).$$
\end{theorem}

An interesting consequence is that the volume  only takes finitely many values on the family of smooth projective varieties of general type with fixed underlying 6-manifold (see Corollary \ref{c_volume}). A second consequence (which follows immediately applying \cite{HM06}) is that the family of all smooth complex projective threefolds of general type with bounded Betti numbers is birationally bounded (see Corollary \ref{c_boundedness}). Such questions remain  open in higher dimension.
In a forthcoming paper, we plan to study the Chern numbers of a variety $Y$ which admits a Mori fibre space structure. 

\medskip

We now describe the second part of our program: we want to determine how the Chern number $c_1^3$ varies under the Minimal Model Program. Recall that if $X$ is a smooth projective threefold and we run a MMP on $X$, then we obtain a birational map $X\dashrightarrow Y$ as a composition of elementary transformations, given by divisorial contractions and flips: 
$$
X=X_0 \dashrightarrow X_1 \dashrightarrow \ldots \dashrightarrow X_m=Y.
$$
We plan to bound $K_{X_k}^3-K_{X_{k-1}}^3$ at each step,  in terms of the topology of the manifold underlying $X$. 

In this paper, we consider the case of divisorial contractions. 
Recall that a divisorial contraction $X_{k-1}\to X_k$ is a birational morphism which contracts a prime divisor $E$ into either a point or a curve. The first  case can be easily handled thanks to Kawakita's classification \cite{Kawhigher}. In particular, we can show that:
$$
0<K_{X_{k-1}}^3-K_{X_k}^3 \le 2^{10} b_2^{2},
$$
where $b_2=b_2(X)$ is the second Betti number of $X$ (see Proposition \ref{p_ctp}).

The case of divisorial contractions to curves is much harder. In general, in this case, the difference between the Chern numbers may not be bounded by a combination of Betti numbers (e.g. consider a blow-up of a rational curve of degree $d$ in $\mathbb P^3$). To deal with this situation we study the integral cubic form $F_{X_i}$ associated to the cup product on $H^2(X_i,\Z)$.
The cubic form $F_X$ is one of the most important topological invariant of a smooth 3-fold $X$ and many topological information of $X$ are encoded in the cubic form $F_{X}$ (e.g. see \cite{OV95}). In the case of a blow-down to a smooth curve $f: W \to Z$ the cubic form $F_W$ assumes a special form
$$
F_W(x_0,\ldots,x_n)= ax_0^3  + 3x_0^2(\sum_{i=1}^n b_i x_i) + F_Z(x_1,\ldots,x_n),
$$ 
which we call \emph{reduced form}. The goal of Section \ref{s_cubics} is to prove a finiteness result on the number of possible reduced forms in the case of cubic forms with non-zero discriminant (see Theorem \ref{t_b2=n}).

In particular,  we can associate to any  projective threefold $X$ a topological invariant $S_X$ which is an integer number depending only on the cubic form $F_X$ of $X$ (see Definition \ref{d_skansen}). 


Our main result is the following. It is obtained by combining together methods in birational geometry, topology and arithmetic geometry. 

\begin{theorem}\label{t_global}
Let $Y$ be a terminal $\Q$-factorial 3-fold with associate cubic form $F_Y$ and let $f: Y \to X$ be a divisorial contraction to a point or to a smooth curve contained in the smooth locus of $X$ (in this last case assume also that $\Delta_{F_Y} \ne 0$). 
\begin{enumerate}
	\item If  $f$ contracts a divisor to a point, then $|K_Y^3-K_X^3| \le  2^{10} b_2(Y)^{2}$. If $f$ contracts a divisor to a curve, then 
	$$
	|K_Y^3-K_X^3| \le 2S_W + 6(b_3(Y)+1),
	$$
 where $S_Y$ is as in definition \ref{d_skansen}. Moreover, the same inequality is true after replacing $b_3(Y)$ by $Ib_3(Y)=\dim IH^3(Y,\Q)$.
	\item The cubic form $F_X$ is determined up to finitely ambiguity by the cubic form $F_Y$.
\end{enumerate}
\end{theorem}

We believe that the methods used to prove Theorem \ref{t_global}  will have interesting applications to questions concerning the topology and the geography of threefolds (see, for example, \cite{BCT16}).

\medskip 

Let $X$ be a smooth threefold and let $f: X \dto Y$ be a minimal model of  $X$. It is very natural to ask which topological invariants of $Y$ are determined by those of $X$. It is known that the Betti numbers of $Y$ are determined up to finite ambiguity by the Betti numbers of $X$ (the case of $b_3$ has been treated very recently in \cite{Chen16}). 

The same question for the ring structure of the cohomology is very delicate. The following immediate consequence of Theorem \ref{t_global} goes in the positive direction.

\begin{corollary}\label{c_top}
Let $X$ be a smooth complex projective threefold. Let $f: X \dto Y$ be a minimal model program for $X$. 

If $f$ is composed only by divisorial contractions to points, then $F_Y$ is determined up to finitely ambiguity by $F_X$.

If $\Delta_{F_X} \ne 0$ and $f$ is a composition of  divisorial contractions to points and blow-downs to smooth curves in smooth loci, then $F_Y$ is determined up to finitely ambiguity by $F_X$.
\end{corollary}

Finally, we can combine the above results to obtain the following corollary.

\begin{corollary}\label{c_main}
Let $X$ be a smooth complex projective threefold which is not uniruled and let $F_X$ be its associated cubic form. Assume that $\Delta_{F_X} \ne 0$ and that there exists a birational morphism $f\colon X\to Y$ onto a minimal projective threefold $Y$, which is  obtained as a composition of  divisorial contractions to points and blow-downs to smooth curves in smooth loci. 

Then there exists a constant $D$ depending only on the topology of the $6$-manifold underlying $X$ such that
$$
|K_X^3| \le D.
$$

\end{corollary}

\medskip

It remains to study divisorial contractions to singular curves and flips. On the other hand, the Minimal Model Program of any smooth projective threefold  may be also factored into a sequence of flops, blow-up along smooth curves and divisorial contractions to points (see \cite{ChenHacon11, Chen13a}). Recall that if $W \dto Z$ is a flop, then $K_W^3=K_Z^3$; thus, it is crucial to study how the cubic form $F$ varies under flops. We will study this problem in a forthcoming paper.

\section{Preliminary Results}
\subsection{Notations}\label{s_notation}
We work over the field of complex numbers. 
We refer to \cite{KM98} for the classical notions  in birational geometry. 
In particular, if $X$ is a normal projective variety, we denote by $K_X$ the {\em canonical divisor} of $X$. We also denote by $\rho(X)$ the \textit{Picard number} of $X$,  by $N^1(X)$ the group of Cartier divisors modulo numerical equivalence and by $\bar {H}^i(X,\mathbb Z)$ the $i$-th singular cohomology group of $X$ modulo its torsion subgroup. In particular, $b_i(X)=\rk {\bar H}^i(X,\mathbb Z)=\dim H^i(X,\mathbb Q)$ is the $i$-th \textit{Betti number} of $X$. We say that $X$ is \textit{$\mathbb Q$-factorial} if every Weil divisor $D$ on $X$ is $\mathbb Q$-Cartier, i.e. there exists a positive integer $m$ such that $mD$ is Cartier. 
If $f\colon Y\to X$ is a birational morphism between normal projective varieties and $K_X$ is $\mathbb Q$-Cartier, then we may write
$$K_Y=f^*K_X+\sum_{i=1}^k a_iE_i$$ 
where the sum is over all the exceptional divisors $E_1,\dots,E_k$ of $f$. The number $a_i$ is the \textit{discrepancy} of $f$ along $E_i$ and it is denoted by $a(E_i,X)$. In particular, $X$ is said to be \textit{terminal} if for any birational morphism $f\colon Y\to X$ and for any exceptional divisor $E$, we have $a(E,X)>0$. Recall that terminal singularities are \textit{rational}, i.e. if $f\colon Y\to X$ is a resolution then $R^if_*\ring Y.=0$ for all $i>0$.  
A terminal variety $X$ is said to be \textit{minimal} if it is $\Q$-factorial and $K_X$ is nef. 

A \textit{contraction} $f\colon Y\to X$ is a proper birational morphism between normal projective varieties. The contraction $f\colon Y\to X$ is said to be \textit{divisorial} if the exceptional locus of $f$ is an irreducible divisor. It is said to be \textit{elementary}, if $\rho(Y)=\rho(X)+1$. Finally, an elementary contraction $f\colon Y\to X$ is said to be $K_Y$-\textit{negative}, if $-K_Y$ is $f$-ample, i.e. the exceptional locus of $f$ is covered by curves $\xi$ such that $K_Y\cdot \xi<0$. Note that if $Y$ is $\mathbb Q$-factorial and $f\colon Y\to X$ is an elementary divisorial contraction, then $X$ is also $\mathbb Q$-factorial. Moreover, if $Y$ is terminal and $f$ is $K_Y$-negative, then $X$ is also terminal.

\begin{definition}\label{d_volume} Let $X$ be a projective variety with terminal singularities. Then, the \textit{volume} of $X$ is given by 
$$\vol(X)=\limsup_{m\to \infty} \frac{n! ~h^0(X,mK_X)}{m^n}
$$\label{def_vol}
where $n$ is the dimension of $X$. 
\end{definition}
In particular, the volume is a birational invariant and if $X$ is a minimal variety of dimension $n$ then 
$$\vol(X)=K_X^n$$
(see \cite[Section 2.2.C]{Lazarsfeld04a} for more details).

\medskip

\subsection{Terminal singularities on threefolds}\label{s_t3} 
We now recall few known facts about terminal singularities in dimension three. 
Let $(X,p)$ be the germ of a three-dimensional terminal singularity. The {\em index} of $p$ is the smallest positive integer $r$ such that $rK_X$ is Cartier. In addition, it follows from the classification of terminal singularities \cite{Mori85}, that 
there exists a deformation of $(X,p)$ into a variety with $h\ge 1$ terminal singularities $p_1,\dots,p_h$ which are isolated cyclic quotient singularities of index $r(p_i)$. The set $\{p_1,\dots,p_h\}$ is called the {\em basket} $\mathcal B(X,p)$ of singularities of $X$ at $p$ \cite{Reid87}. 
As in \cite{ChenHacon11}, we define 
$$\Xi(X,p)=\sum_{i=i}^h  r(p_i).$$
Thus, if $X$ is a projective variety of dimension $3$ with terminal singularities and $\Sing X$ denotes the finite set of singular points of $X$, we may define
$$\Xi(X)=\sum_{p\in\Sing X} \Xi(X,p).$$

\begin{lemma}\label{l_xi}
Let $(X,p)$ be the germ of a three-dimensional terminal singularity and let $\mathcal B(X,p)$ be the basket at $p$. 

Then, for each $q\in \mathcal B(X,p)$,  the index $r(q)$ of $q$ divides $4\cdot \Xi(X,p)$. 
\end{lemma}
\begin{proof}
It follows from the classification of terminal singularities, that the points of the basket $\mathcal B(X,p)$ either have all the same index $r$ or their index divides $4$ when $r(p)=4$ and $p \in X$ is of type $cAx/4$ (e.g. see \cite[Remark 2.1]{ChenHacon11}). Thus the claim follows. 
\end{proof}

By \cite[Proposition 3.3]{CZ14}, we have:
\begin{proposition}\label{p_2b2}
Let $X$ be a smooth projective threefold and 
assume that 
$$X=\rmap X_0.\rmap .\dots. .X_k=Y$$ 
is a sequence of steps for the $K_X$-minimal model program of $X$. 

Then 
$$\Xi(Y)\le 2b_2(X).$$
In particular, the inequality holds if $Y$ is the minimal model of $X$. 
\end{proposition}

\medskip

In the proof of our main results, we will use the Bogomolov-Miyaoka-Yau inequality and the Riemann Roch formula for terminal threefolds. Recall that, on any terminal threefold $X$,  we  may define  $c_1(X)$ as the anti-canonical divisor $-K_X$ and, for any $\Q$-Cartier divisor $D$ on $X$ we define the number $D.c_2(X)$ as $f^*D.c_2(Y)$ where $f\colon Y \to X$ is any resolution of $X$. It is easy to check that the definition does not depend on the resolution. 

\begin{theorem}\label{t_BMY}
Let $Y$ be a minimal three-dimensional projective variety with terminal singularities. 

Then 
$$
(3c_2 - c_1^2).c_1 \le 0.
$$
\end{theorem}
\begin{proof}
It follows from  \cite[Theorem 1.1]{Miyaoka87b}.
\end{proof}

\begin{theorem}\label{t_RR}
Let $Y$ be a three-dimensional projective variety with terminal singularities. 

Then the holomorphic Euler characteristic of $Y$ is given by 
$$\chi(Y,\ring Y.)= \frac {1} {24} (-K_Y\cdot c_2(Y)+ e)$$
where 
$$e=\sum_{p_\alpha} \left (r(p_\alpha)-\frac 1 {r(p_\alpha)}\right ),$$
and the sum runs over all the points of all the baskets of $Y$. 
\end{theorem}

\begin{proof} See
 \cite{Kawamata86,Reid87}.
\end{proof}

\medskip

\subsection{Cubic Forms}\label{s_cf}

For any polynomial $P\in \mathbb C[x_0,\dots,x_n]$, we denote by $\partial_i P(x)$ the partial derivative of $P$ with respect to $x_i$ at the point $x\in \mathbb C^{n+1}$.  For any ring $R\subseteq \mathbb C$ and for any positive integer $d$, we denote by $R[x_0\dots,x_n]_d$ the set homogeneous polynomials of degree $d$ with coefficients in $R$. 

Given a cubic form $F\in \mathbb C[x_0,\dots,x_n]$, i.e. an homogeneous polynomial of degree $3$, let 
$$\mathcal H_F(x)=(\partial_i\partial_jF(x))_{i,j}$$ be the Hessian of $F$ at the point $x\in \mathbb C^{n+1}$. Note that, for any $x\in \mathbb C^{n+1}$ and for any $\lambda\neq 0$, the rank of $\mathcal H_F$ at the point $\lambda x$ is constant with respect to $\lambda$ and therefore we will denote, by abuse of notation, $\rk \mathcal H_F(p)$ to be the rank of $\mathcal H_F$ at any point in the class of $p \in \mathbb P^n$. 
We say that $F$ is {\em non-degenerate} if $\rk \mathcal H_F$ is maximal at the general point of $\mathbb P^n$, i.e. if $\det \mathcal H_F$ is not identically zero.

\medskip
   
 Let $F(x_0,\ldots,x_n)=\sum_I c_{I}x^{I}\in \mathbb C[x_0,\dots,x_n]_d$. Then the \textit{discriminant} $\Delta_F$ of $F$ is the unique (up to sign) polynomial with integral  coefficients in the  variables $c_I$ such that $\Delta_F$ is irreducible over $\Z$ and  $\Delta_F=0$ if and only if the hypersurface $\{F=0\} \subseteq \mathbb P_{\C}^n$ is singular (see \cite[pag. 433]{GKZ} for more details). In particular, the discriminant  is an invariant under the natural $\SL(n+1,\mathbb C)$-action.

If $F \in \mathbb C[x,y,z]$ is a ternary cubic form, then we denote by $S_F$ and $T_F$ the two $\SL(3,\mathbb C)$-invariants of $F$ as defined in  \cite[ 4.4.7 and 4.5.3]{Sturmfels93}. Then the discriminant of $F$ satisfies
$$
\Delta_F=T_F^2-64S_F^3.
$$

\begin{lemma}\label{l_disc}
Let $F\in \Z[x_0,\dots,x_n]_3$ be an integral cubic form and assume that
$$
F(x_0,\ldots,x_n)=ax_0^3 + x_0^2(\sum_{i=1}^{n}b_ix_i)+ G(x_1,\ldots,x_n)
$$
for some $G\in \Z[x_1,\dots,x_n]_3$. Then $\Delta_G$ divides  $\Delta_F$.
\end{lemma}

\begin{proof}
If $P$ is a polynomial with integral coefficients we denote by $\mathrm{ct}(P)$ the \textit{content} of $P$, that is the gcd of the coefficients of $P$. As in the case of one variable, it is  easy to see that the content is multiplicative.

Let $A$, $\{B_i\}_{i=1,\ldots,n}$ and $\{C_J\}$ be variables and consider the cubic form
$$
f=Ax_0^3 + x_0^2(\sum_{i=1}^{n}B_ix_i)+ g(x_1,\ldots,x_n)
$$
where $g=\sum_J C_{J}x^{J}$. Then $\Delta_f$ and $\Delta_g$ are polynomial in $\Z[A, B_i, C_J]$.  We want to show that $\Delta_g$ divides $\Delta_f$.

Let $R=\C[A,B_i,C_J]$ and let $Z(f), Z(g) \subseteq \mathbb P^N_{\C}=\Proj R$ be the closed subsets defined by $\Delta_f=0$ and $\Delta_g=0$ respectively. Note that $Z(g) \subseteq Z(f)$ because if $\{g=0\}$ has a singular point $z=[z_1,\ldots,z_n]$, then $[0,z_1,\ldots,z_n]$ is a singular point of $\{f=0\}$.  Since $\Delta_g$ is irreducible over $\mathbb Q$ by definition, and hence $Z(g)$ is reduced over $\C$, we deduce that $\Delta_f=\Delta_g \cdot H$ where $H \in R$. 

We need to show that  $H \in \ZZ[A,B_i,C_J]$. We proceed as in the proof of Gauss lemma. We start assuming by contradiction that $H \notin \mathbb Q[A,B_i,C_J]$. Fix an order on $R$ and consider the maximal monomial $m$ in $H$ such that its coefficient is not rational. Consider now the product between $m$ and the highest monomial in $\Delta_g$ to get a contradiction.   Hence $H \in \mathbb Q[A,B_i,C_J]$. 

The claim follows from  the fact that the content of $\Delta_g$ is 1 and that the content is multiplicative.
\end{proof}

\medskip

We have:
\begin{lemma}\label{l_trivialh}
Let $F\in \mathbb C[x_0,\dots,x_n]$ be a cubic form such that there exists a point $p\in \mathbb P^n$ for which 
$\rk \mathcal H_F(p)=0$ (i.e. $\mathcal H_F(p)$ is the trivial matrix).

Then after a suitable coordinate change, $F$ depends on at most $n$ variables. In particular, $\det \mathcal H_F$ vanishes identically on  $\mathbb P^n$. 
 \end{lemma}
\begin{proof}
Euler's formula for homogeneous polynomials implies  that $$F(p)=\partial_{i}F(p)=0\quad\text{ for all }i=0,\dots,n.$$ After a suitable coordinate change, we may assume that $p=(1,0,\dots,0)$. Let $f(y_1,\dots,y_n)=F(1,y_1,\dots,y_n)$. By Taylor's formula,  $f$ is a homogeneous polynomial of degree $3$. Thus, $F(x_0,\dots,x_n)=f(x_1,\dots,x_n)$ and the claim follows. 
\end{proof}

As mentioned in the introduction, arithmetic geometry will play an important role for the proof of our main theorem. In particular, we need the following: 
\begin{theorem}[Siegel Theorem]\label{t_siegel}
Let $R$ be a ring finitely generated over $\mathbb Z$. Let $C$ be an affine smooth curve defined over $R$ and of  genus  $g \ge 1$. 

Then there are only finitely many $R$-integral points on $C$. 
\end{theorem}
\begin{proof}
See \cite[Ch. 8, Theorem 2.4]{Lang83}.
\end{proof}

\medskip

\subsection{Reduced triples}\label{s_reduced} Given a ring $A$, we denote by $\mathcal M(n, A)$ the set of all matrixes with coefficients in $A$, by $\GL(n,A)$ the subgroup of invertible matrixes and by $\SL(n,A)$ the subgroup of matrixes with determinant 1. 

Given a cubic form $F \in \mathbb C[x_0,\dots,x_n]$ and a matrix $T\in \GL(n+1,\mathbb C)$, we will denote by $T\cdot F$ the cubic form given by 
$$T\cdot F(x) = F( T\cdot x).$$

We define
$$W_F=\{ p \in \mathbb P^n\mid \rk \mathcal H_F(p)\le 1\}$$
and 
$$V_F=\{p\in \mathbb P^n\mid \rk \mathcal H_F(p)\le 2\}.$$

\begin{definition}\label{d_reduced}
Let $F\in  R[x_0,\dots,x_n]$ be a non-degenerate cubic form where $R$ is a commutative ring.  We say that $(a,B,G)$ is a {\em reduced triple} associated to $F$  if there exists an element $T \in \SL(n+1,R)$ such that
\begin{equation}\label{e_rf}
T \cdot F=ax_0^3 + x_0^2\cdot \sum_{i=1}^n b_ix_i +G(x_1,\dots,x_n)
\end{equation}
where $a \in R$, $B=(b_1,\ldots,b_n) \in R^n$ and $G\in R[x_1,\dots,x_n]$ is a non-degenerate cubic form. 
For simplicity, we will denote \eqref{e_rf} as 
$$T \cdot F=(a,B,G).$$ 
 In this case we also say that $T\cdot F$ \emph{is in reduced form} $(a,B,G)$.

We say that two reduced triples $(a,B,G)$ and $(a',B',G')$, are {\em equivalent over $R$} if $a=a'$ and there is an element $M \in \SL(n, R)$ such that  $B'=M\cdot B$ and $G'=M \cdot G$. 
\end{definition}

The motivation to study the loci $W_F$ and $V_F$ and reduced forms comes from Propositions \ref{p_cubicpoints} and \ref{l_F_Y}. More precisely, it is easy to see that if  $F\in \mathbb C[x_0,\dots,x_n]$ is a cubic form in reduced form, and  $p=[1,0,\dots,0]$, then $p\in V_F$ (see for example \cite[Lemma 2.1]{BCT16}).

In the subsequent we will use the following result:

\begin{theorem}[Jordan's theorem] \label{t_jordan}
Let $F\in \mathbb Z[x_0,\dots,x_n]_3$ be a cubic form with non-zero discriminant $\Delta_F$ and consider the set 
$$\mathcal A_F = \{T\cdot F \mid T\in \SL(n+1,\C) \}\subseteq \C[x_0,\dots,x_n]_3.$$

Then the quotient 
$$(\mathcal A_F\cap \mathbb Z[x_0,\dots,x_n]_3)/ \SL(n+1,\mathbb Z)$$
is finite. 
%
\end{theorem}
\begin{proof}
It follows from \cite[Corollary 4 and 5]{OV95}.
\end{proof}

\medskip

\subsection{Cubic forms on threefolds}\label{s_cft}

Let $X$ be a terminal $\Q$-factorial projective threefold. 
Let $\underline{h}=(h_1, \ldots, h_n)$ be a basis of $\bar{H}^2(X, \Z)$. The intersection cup product induces  a symmetric trilinear form 

$$
\phi_X: \bar{H}^2(X,\Z) \otimes \bar{H}^2(X,\Z) \otimes \bar{H}^2(X,\Z) \to H^6(X,\Z) \cong \Z.
$$
Thus, we may define a cubic homogeneous polynomial $F_X \in \Z[x_1,\ldots,x_n]$ as
$$
F_X(x)= \sum_{\substack{I=(i_1,\ldots,i_n): \\  i_1+\ldots +i_n=3}} \binom{3}{I}  \phi_X(\underline{h}^{I}) x^{I}.
$$  
We call $F_X$ \emph{the cubic form associated to $X$}. 

As in the smooth case, we have:

\begin{lemma}
The cubic form $F_X$ is non-degenerate, that is  $\det \mathcal H_{F_X}$ is not identically zero. 
\end{lemma}

\begin{proof}
Let $\Sigma \subseteq X$ be the singular locus of $X$. Since $X$ is terminal, $\Sigma$ is a finite set and there exists  a resolution $\pi: Y \to X$ with divisorial exceptional locus $E$  such that $Y \setminus E$ is isomorphic to $X \setminus \Sigma$. 

Let $\{\gamma_0, \ldots, \gamma_b\}$ be a basis of  $H^2(X,\Q)$ and let $\mathcal B=\{\beta_i=f^*\gamma_i\}$. After completing $\mathcal B$ to a basis  of $H^2(Y,\Q)$, we may write
$$
F_Y(x_0, \ldots, x_n)= F_X(x_0, \ldots, x_b) + F(x_{b+1}, \ldots, x_n),  
$$ 
where we are considering the cubic forms over $\Q$.

\cite[Proposition 16]{OV95} implies that $\det \mathcal H_{F_Y}$ is not identically zero. Since $\det \mathcal H_{F_Y}=\det \mathcal H_{F_X } \cdot \det \mathcal H_{F}$, the claim follows. 
\end{proof}

\begin{definition}\label{d_skansen}
Let $X$ be a terminal $\mathbb Q$-factorial projective threefold and let $F_X\in \mathbb Z[x_1,\dots,x_n]_3$ be the cubic form associated to $X$. We define
$$
S_X:=\sup\{ |a|\in \mathbb Z \mid \mbox{there exists } T \in SL(n+1,\Z) \mbox{ s.t. } T\cdot F_X=(a,B,G)  \},
$$
where we set $S_X=0$ if there are no reduced triples associated to  $F_X$. 
\end{definition}

Note that $S_X$ is a  topological invariant of $X$ since $F_X$ is a topological invariant (modulo the action of $\SL(n+1,\Z))$. 

\medskip

\subsection{Topology of threefolds}

We now study how the Betti numbers behave under a birational morphism (see \cite{Caibar05} for some related results).
Note that the singularities of a  $\Q$-factorial terminal threefold $X$ are in  general not analytically $\Q$-factorial. In particular, $X$ is  in general not a $\Q$-homology manifold (see \cite[Lemma 4.2]{Kollar89}) and  the singular cohomology may differ from  the intersection cohomology. 

In dimension three,  all the  Betti numbers behave well under birational transformations except for $b_3$ (see Lemma \ref{l_betti}). The behaviour of the third Betti number is more subtle and depends on the singularities of $X$ and $Y$ as the following example shows:

\begin{example}
Let $X\subseteq \mathbb P^4$ be a quartic  threefold with just one node (rational double point) $p\in X$. It is known that $X$ is $\Q$-factorial (e.g. see  \cite{Cheltsov06}). Locally, the germ $(X,p)$  may be written as 
$$
\{xy-wz=0\} \subseteq \C^4,
$$
which is not analytically $\Q$-factorial.
Let $f\colon Y \to X$ be the blow-up of the singularity and let $E \cong  \mathbb P^1 \times \mathbb P^1$ be the exceptional divisor. It follows that
$$
b_3(Y)=b_3(X)-1.
$$
\end{example} 

\medskip 

In particular, the third Betti number may increase under some of the steps of the Minimal Model Program. For this reason, it will be often useful to look at the intersection cohomology instead. 

Given a projective variety $X$, we denote by $IH^i(X,\Q)$ the \textit{middle-perversity intersection cohomology group} of dimension $i$ and by $Ib_i$ its dimension. Note that if $X$ is smooth then $IH^i(X,\Q)$  coincides with $H^i(X,\Q)$ and in particular $Ib_i(X)=b_i(X)$ for all $i$. 

 We will use the following  consequence of the decomposition theorem for intersection cohomology (see \cite{BBD}):

\begin{theorem} \label{t_decomposition}
Let $f\colon Y \to X$ be a proper birational morphism between algebraic varieties. Assume that $Y$ is smooth. Then the cohomology $H^*(Y,\Q)=IH^*(Y,\Q)$ of $Y$ contains the intersection cohomology $IH^*(X,\Q)$ of $X$ as a direct summand.
\end{theorem}

We now restrict our study to  the case of threefolds:

\begin{lemma}\label{l_cohomology} 
Let $f\colon Y\to X$ be a birational  morphism between  projective threefolds with terminal singularities. Let $E$ be an exceptional divisor of $f$ and let $W=f(E)$.  Assume that  $f$ induces an isomorphism $Y \setminus E \to X \setminus W$. 

Then
$$
0 \to H^i(X,\Q) \to H^i(Y,\Q) \oplus H^i(W,\Q) \to H^i(E,\Q) \to 0 
$$
is exact for any $i \ge 4$  and
$$
0 \to IH^i(X,\Q) \to IH^i(Y,\Q) \oplus IH^i(W,\Q) \to IH^i(E,\Q) \to 0 
$$
is exact for any $i \ge 1$.

\end{lemma}
\begin{proof}
  From the exact sequence of the pairs we get a long exact sequence in  cohomology
$$\cdots \to H^i(X,\Q) \to H^i(Y,\Q) \oplus H^i(W,\Q) \to H^i(E,\Q) \to H^{i+1}(X,\Q)\to \cdots
$$which by \cite[Prop. 8.3.9]{Deligne74} is an exact sequence of mixed Hodge structure.

Since $X,Y$ have isolated singularities, for $i \ge 4$ the Hodge structure on $H^i(X,\Q)$  is pure of weight $i$ (see \cite{Steenbrink83}). On the other hand, since $E$ is projective,  $H^k(E,\Q)$ has weight at most  $k$ for any $k$ (\cite[Thm. 8.2.4]{Deligne74}). Thus,  the maps
$$
H^{i}(E,\Q) \to H^{i+1}(X,\Q)
$$
are zero for $i\ge 3$. 

The same argument applies for intersection cohomology with the advantage that  the Hodge structure on $IH^i(X,\Q)$ is pure of weight $i$ for any $i$ by \cite{Saito88}.
\end{proof}

\begin{lemma}\label{l_betti} Let $f\colon Y \to X$ be an elementary divisorial  contraction between  $\mathbb Q$-factorial projective threefolds with terminal  singularities.

Then
\begin{enumerate}
\item $b_0(Y)=b_6(Y)=b_0(X)=b_6(Y)=1, $
\item $ b_1(Y)=b_1(X), $
\item $b_2(Y)=b_2(X) + 1 $
\item $ b_4(Y)=b_4(X) + 1$, and 
\item $b_5(Y)=b_5(X).$ 
\end{enumerate}
\end{lemma}
\begin{proof}
$(1)$ is clear. Lemma \ref{l_cohomology} implies $(4)$ and $(5)$.

\smallskip

%
%

We now want to show that $R^1f_*\Z=0$. It is enough to show it locally around any point $x\in X$. 
We consider the exact sequence
$$\begin{aligned}
0\to f_*\Z \to & f_*\ring Y.\stackrel {\exp}\longrightarrow  f_*\mathcal O^*_Y\\
\to R^1f_*\Z&\to R^1f_*\ring Y.
\end{aligned}$$
The exponential map is surjective locally around $x\in X$.  Since $X$ and $Y$ have rational  singularities, it follows that $R^1f_*\ring Y.=0$. Thus, $R^1f_*\mathbb Z=0$, as claimed. The Leray spectral sequence implies that $H^1(X,\Z)\to H^1(Y,\Z)$ is an isomorphism and, in particular, $(2)$ follows. 


Let $H_2(Y/X,\mathbb C)\subseteq H_2(Y,\mathbb C)$ be the subspace generated by all the images of $H_2(F,\mathbb C)$, where $F$ runs through all the fibres of $f$. 
 \cite[Theorem 12.1.3]{KM92} implies  that  $H_2(Y/X,\C)$ is generated by algebraic cycles and that there exists an exact sequence:
$$0\to H_2(Y/X,\mathbb C)\to H_2(Y,\mathbb C)\to H_2(X,\mathbb C)\to 0.$$
Since $f$ is an elementary divisorial contaction, it follows that  all the non-trivial algebraic cycles contained in the fiber of $f$ are numerically proportional to each other and, in particular,   $$\dim H_2(Y/X,\C)=1.$$ Thus, $(3)$ follows. 
\end{proof}

\section{Cubic forms in reduced form}\label{s_cubics}

The aim of this section is to prove the following:

\begin{theorem}\label{t_b2=n}
Let $F\in \ZZ[x_0,\dots,x_n]$ be a non-degenerate cubic form (cf. \S \ref{s_cf}) with non-zero discriminant $\Delta_F$. 

Then there are  finitely many triples 
$$(a_i,B_i,G_i)\in \mathbb Z\times \mathbb Z^n\times \mathbb Z[x_1,\dots,x_n]_3\qquad i=1,\dots,k$$
such that any reduced triple associated to $F$ is equivalent to 
$(a_i,B_i,G_i)$ over $\mathbb Z$ for some $i \in \{1,\dots,k\}$ (cf. Definition \ref{d_reduced}). 

In addition, we have that $\Delta_{G_i} \ne 0$ for all $i=1,\dots,k$.
\end{theorem}

Before we proceed with the proof of Theorem \ref{t_b2=n}, we first sketch some of its main ideas. Note 
that if $F$ is in reduced form $(a,B,G)$ then the point $p=(1,0,\dots,0)$ is contained in the set $V_F$, 
defined in \S \ref{s_reduced} . Thus, our first goal is to show that the set of points  $p\in V_F$ such that 
$F(p)\neq 0$ is contained in a finite union of points, lines and plane cubics (cf. Theorem \ref{t_plane}). 
Assuming furthermore that the discriminant $\Delta_F$ of $F$ is not zero, we characterise the cubic forms $F$ which contain a line (cf. Corollary \ref{c_Cline}) or a plane curves 
(cf. Corollary \ref{c_Ccubic}) inside $V_F$.

The next step is to restrict the cubic form to one of the lines or 
plane curve contained in $V_F$. To deal with this situation we study binary (cf. Proposition \ref{p_binary}) and ternary 
cubic forms (cf. Proposition \ref{p_b2=3}) with non-zero discriminant. The main tool used in the proof of 
these results is Siegel's  theorem on the finiteness of integral  points in a curve 
of positive genus. Finally, we conclude the proof of Theorem \ref{t_b2=n} in \S \ref{s_general}. 

\medskip

\subsection{Points of low rank for a cubic form}
In this subsection, we  study the sets $W_F$ and $V_F$ associated to a cubic form $F\in \mathbb C[x_0,\dots,x_n]$ (cf. \S \ref{s_reduced}). Many of the result below  depend on some simple calculations on  cubics forms. To illustrate some of the methods presented below, we begin  with a basic result:
\begin{lemma}\label{l_easy}
Let 
$$F=x^3_0+x_0Q+R\in \mathbb C[x_0,\dots,x_n]_3$$
 be a  cubic form, where $Q, R\in \mathbb C[x_1,\dots,x_n]$ are homogeneous polynomials of degree $2$ and $3$ respectively. Let $A$ be the $n\times n$ symmetric matrix associated to $Q$. Let $p=[1,0,\dots,0]$. 

Then $\rk \mathcal H_F(p)=\rk A+1$. 
\end{lemma}
\begin{proof}
The claim is a simple computation. 
\end{proof}

We now proceed by studying the set $W_F$ (cf. \S \ref{s_reduced}) associated to a non-degenerate cubic form $F$:

\begin{proposition}\label{p_cubic}  
Let $F\in \mathbb C[x_0,\dots,x_n]$ be a non-degenerate cubic form. 
Then $W_F$ is a finite set. 
\end{proposition}

\begin{proof}
Let $W'_F=W_F\cap \{F=0\}$.
We first show that $W'_F$ is a finite set. 
Assume by contradiction that there exist an irreducible curve $C$ inside $W'_F$ and 
let $p\in C$.
We  say that an hyperplane $H\subseteq \mathbb P^n$ is {\em associated} to $p$ if:
\begin{enumerate}
\item $\det \mathcal H_F$ vanishes along $H$,
\item $p\in H$, and 
\item if $G=F_{|H}$ then  $\mathcal H_G(p)$ is trivial.
\end{enumerate}
Lemma \ref{l_trivialh} implies that  $\rk \mathcal H_F(p)=1$. 
After taking a suitable coordinate change, we may assume that 
$p=[1,0,\dots,0]$. In particular 
$$F(x_0,\dots,x_n)= x_{0}^2\cdot L_1+ x_0 \cdot Q_1+ R_1 $$
for some homogeneous polynomials $L_1,Q_1,R_1\in\mathbb C[x_1,\dots,x_{n}]$ of degree $1,2$ and $3$ respectively. Since $p\in W_F$, it follows that $L_1=0$. By assumption, $Q_1$ is not zero. Using again the fact that $p\in W_F$, similarly to Lemma \ref{l_easy}, it follows that, after taking a suitable coordinate change in $x_1,\dots,x_{n}$, we may assume that $Q_1=x_{1}^2$. 
We may write 
$$R_1(x_1,\dots,x_{n})= x_{1}^2\cdot L + x_{1} \cdot Q + R $$
for some homogeneous polynomials $L\in\mathbb C[x_1,\dots,x_{n}]$ and $Q,R\in\mathbb C[x_2,\dots,x_{n}]$ of degree $1,2$ and $3$ respectively.
After replacing $x_{0}$ by $x_0+L$, we may assume that $L=0$. 
Thus, we have
$$F(x_0,\dots,x_n)=x_0\cdot x_{1}^2 + x_{1}\cdot Q + R.$$
Let $H_p=\{x_{1}=0\}$.  An easy computation shows that $H_p$ is an hyperplane associated to $p$. We now show that  
such an hyperplane is unique. Assume that $H'\subseteq \mathbb P^{n}$ is also an hyperplane associated to $p$. Since $p\in H'$, we have $H'=\{\ell=0\}$ for some linear function $\ell\in \mathbb C[x_1,\dots,x_{n}]$. If $H'\neq H_p$, after a suitable change of coordinates in $x_2,\dots,x_{n}$, we may assume that 
$$\ell=x_n-\alpha x_{1}$$
for some $\alpha \in \mathbb C$. 
Thus if $G'=F_{|H'}$, we may write
$$G'(x_0,\dots,x_{n-1})=x_0x_{1}^2  + x_{1}  Q(x_2,\dots,x_{n-1},\alpha x_{1})+R(x_2,\dots,x_{n-1},\alpha x_{1})$$
and it follows that 
$$\partial_{1}\partial_{1} G'(p)\neq 0$$
which contradicts (3). Thus, $H'=H_p$ and the claim follows. 

Now let $q\in C$ be a point such that $H_p=H_q$. We want to show that  $q=p$. If $R=0$ then it follows easily that $W'_F=\{p\}$. Thus, by Lemma \ref{l_trivialh}, after a suitable change in coordinates in $x_2,\dots,x_{n}$, we may assume that $R=R(x_{n-k},\dots,x_{n})$ for some $k\ge 0$ and that there is no point $z\in \mathbb P^k$ such that $\mathcal H_R(z)$ is trivial. If $q=[y_0,\dots,y_{n}]$, it follows by (3) that 
$$y_{n-k}=\dots=y_{n}=0.$$ 
Since $\rk \mathcal H_F(q)=1$, it follows the that the minor spanned by the $i$-th and $(n-i)$-th rows and columns of $\mathcal H_F(p)$ must have trivial determinant for any $i=0,\dots,n-2$ and in particular, since $y_{1}=0$ and $\mathcal H_R(y_2,\dots,y_{n})$ is trivial, 
it follows that $\partial_iQ(y_0,\dots,y_n)=0$. It is easy to show that 
this implies that if $q\neq p$ then $\det \mathcal H_R$ vanishes identically, a contradiction. 

Since by assumption $\det \mathcal H_F$ is a non-trivial function, there exist only finitely many hyperplanes on which  $\det\mathcal H_F$ vanishes and  (1) implies that $H_p=H_q$ for infinitely many pair of points $p, q\in C$, a contradiction. 
Thus, $W'_F$ is a finite set. 

Now let $p\in W_F$ be a point such that $F(p)\neq 0$. After a suitable change of coordinates, we may assume that $p=[1,0,\dots,0]$ and that 
$$F(x_0,\dots,x_n)=x_0^3 + x_0^2\cdot L + x_0\cdot Q + R$$
for some homogeneous polynomials  $L,Q,R\in \mathbb C[x_1,\dots,x_{n}]$ of degree $1$, $2$ and $3$ respectively. 
After replacing $x_0$ by $x_0+\frac 1 3 L$ we may assume that $L=0$. Since $p\in W_F$, Lemma \ref{l_easy} implies  that $Q=0$.
Let $q=[z_0,\dots,z_n]\in W_F$.  Then either $q=p$ or $z_0=0$ and $[z_1,\dots,z_{n}]\in W_R$. Thus, the result follows by induction on $n$. 
\end{proof}

\begin{remark}
Note that the same result does not hold if we replace the assumption that $F$ is non-degenerate, by  the weaker assumption that $\rk \mathcal H_F(p)\ge 1$ for any $p\in \mathbb P^n$ (see Lemma \ref{l_trivialh}). E.g. consider $$F(x_0,\dots,x_4)=x_4x_3^2+x_3x_1x_0+x_2x_1^2.$$
Then it is easy to check that $W_F$ is not finite. 
\end{remark}

\medskip

We now proceed by studying the set $V_F$ (cf. \S \ref{s_reduced}) associated to a non-degenerate cubic form  $F\in \mathbb C[x_0,\dots,x_n]$. More specifically, if $V_F$ contains a curve $C$ on which $F$ is not identically zero, then we may write $F$ in a normalised  form as in Theorem \ref{t_special}. The result will be crucial in our proof of  Theorem \ref{t_plane} below. In order to obtain a normalisation as in Theorem \ref{t_special}, we proceed similarly as in the proof of Proposition \ref{p_cubic}. Indeed, by Lemma \ref{l_easy}, to any  point $p\in C$ such that $F(p)\neq 0$, we may associate an hyperplane in $\mathbb P^n$ which contains $p$. The normalisation of $F$ will then depend on whether the curve $C$ is contained in this hyperplane or not. 

Fix a positive integer $n$ and let $\ell$ and $k$ be non-negative integers such that $n\ge \ell+2k+1$. We will denote:
$$I_{\ell,k}=\{\ell+2i+1\mid i=0,\dots,k\}\cup\{\ell+2k+2,\dots,n \}.$$
Given a finite subset $I\subseteq \mathbb N$, we will also denote by
$\mathbb C[x_I]$ the algebra of polynomials in $x_i$ with $i\in I$.

\begin{theorem}\label{t_special}
Let $F\in \mathbb C[x_0,\dots,x_n]$ be a non-degerate cubic form. Let $C\subseteq V_F$ be a curve such that $F(p)\neq 0$ at the general point of $C$. 

Then, there exist non-negative integers $\ell,k$ such that,   after a suitable change of coordinates, we may write
$$F=\sum_{i=0}^{\ell} G_i +\sum_{i=1}^k (x^2_{\ell+2i+1} + M_i)\cdot x_{\ell + 2i} + R_{\ell+k+1}$$
where 
\begin{enumerate}
\item\label{i_G} $G_i\in \mathbb C[x_i,x_{i+1}]$ is a cubic form for any $i=0,\dots,\ell$ with $$G_0=x_0^3+x_0x_1^2;$$
\item\label{i_M} $M_i= \delta_i x^2_{\ell+1}$ for any $i=1,\ldots,k$ with $\delta_i \in \mathbb C$;
\item \label{i_R} $R_{\ell+k+1}\in \mathbb C[x_{I_{\ell,k}}]$ is a cubic form;
\item $C\subseteq\bigcap_{i\in I_{\ell,k+1}}\{x_i=0 \}$. 
\end{enumerate}

Moreover if $C \not\subseteq \{ x_{l+2k+2}=0 \}$  we may write
$$
R_{\ell+k+1}= M_{k+1}\cdot x_{\ell+2k+2}+R_{l+k+2}
$$ 
where
\begin{enumerate}
\setcounter{enumi}{4}
\item \label{i_k} $R_{\ell+k+2}\in \mathbb C[x_{I_{\ell,k+1}}]$ is a cubic form and $M_{k+1}\in\mathbb C[x_{\ell+1},x_{\ell+3},\dots,x_{\ell+2k+1}]$ is a quadric.
\end{enumerate}

\end{theorem}

\begin{proof} We divide the proof in $4$ steps: 

\medskip 

\textbf{Step 1.} 
By Proposition \ref{p_cubic} there exists 
$p\in C$ such that $F(p)\neq 0$ and $\rk \mathcal H_F(p)=2$.
Since $F(p)\neq 0$, after a suitable change of coordinates we may assume that
$p=[1,0,\dots,0]$ and 
$$F=x_0^3 + x^2_0L + x_0Q + R$$
for some homogeneous polynomials $L,Q,R\in \mathbb C[x_1,\dots,x_{n}]$ of degree $1$,$2$ and $3$ respectively.
After replacing $x_0$ by $x_0-\frac 1 3 L$ we may assume that $L=0$. Since $\rk\mathcal H_F(p)=2$, by Lemma \ref{l_easy},  after a suitable change of coordinates in $x_1,\dots,x_{n}$, we may assume that $Q=x_1^2$. Thus, we have
$$F=G_0 + R_1,$$
where $G_0=x_0^3+x_0 x_1^2$ and $R_1=R\in \mathbb C[x_1,\dots,x_n]$. We distinguish two cases. 
If $C$ is contained in the hyperplane $\{x_1=0\}$, then we set $k=\ell=0$ and we continue to Step 3. Otherwise, we set $\ell=1$ and we proceed to Step 2. 

\medskip 

\textbf{Step 2.} 
We are assuming that 
$$F=\sum_{i=0}^{\ell-1} G_i + R_\ell$$
where $G_i\in \mathbb C[x_i,x_{i+1}]$ and $R_\ell\in \mathbb C[x_{\ell},\dots,x_n]$ are cubic forms, 
and $C$ is not contained in the hyperplane $\{x_{\ell}=0\}$. 
We claim that after a suitable change of coordinates in $x_{\ell},\dots,x_n$, we may write 
$$R_\ell=G_{\ell}+R_{\ell+1}$$
where $G_{\ell}\in \mathbb C[x_{\ell},x_{\ell+1}]$ and $R_{\ell+1}\in  \mathbb C[x_{\ell+1},\dots,x_n]$ are cubic forms.
Assuming the claim, if $C$ is contained in the hyperplane $\{x_{\ell+1}=0\}$ we set $k=0$ and we proceed to Step 3. Otherwise, we replace $\ell$ by $\ell+1$ and we repeat Step 2. 

We now prove the claim. By assumption, there exists $q\in C$ such that $q\notin \{x_{\ell}=0\}$.
After a suitable change of coordinates in $x_{\ell},\dots,x_n$, we may assume that 
$$q=[z_0,\dots,z_{\ell-1},1,0,\dots,0],$$ 
for some $z_0,\dots,z_{\ell-1}\in \mathbb C$. 
We may write
$$R_{\ell}=\alpha_\ell x_\ell^3+L_\ell x_\ell^2+Q_\ell x_\ell+R_{\ell+1},$$
for some  homogeneous polynomials $L_\ell,Q_\ell,R_\ell\in \mathbb C[x_{\ell+1},\dots,x_{n}]$ of degree $1$,$2$ and $3$ respectively.
Since $\rk \mathcal H_F(q)\le 2$, 
after a suitable change of coordinates, we may write  $L_\ell=\beta_\ell x_{\ell+1}$ and $Q_\ell=\gamma_\ell x^2_{\ell+1}$
for some $\beta_\ell,\gamma_\ell\in \mathbb C$. We may define 
$$G_{\ell}=\alpha_\ell x_\ell^3+\beta_\ell x_\ell^2 \cdot x_{\ell+1} + \gamma_\ell x_\ell\cdot x_{\ell+1}^2$$
and the claim follows. 

\medskip 

\textbf{Step 3.}
We are assuming that
$$F=\sum_{i=0}^{\ell} G_i +\sum_{i=1}^k (x^2_{\ell+2i+1} + M_i)\cdot x_{\ell + 2i} + R_{\ell+k+1}$$
where $G_i,$ $M_i$ and $R_{\ell+k+1}$ satisfy \eqref{i_G}, \eqref{i_M} and \eqref{i_R} and $$C\subseteq \{x_{\ell+1}=x_{\ell+3}=\dots=x_{\ell+2k+1}=0\}.$$
If we also have that 
$$C\subseteq \{x_{\ell+2k+2}=\dots=x_n=0\}$$
then we are done. 
In particular, if $n<\ell+2k+2$, then we are done. 
Otherwise, after a suitable change of coordinates in $x_{\ell+2k+2},\dots,x_{n}$ we may assume that 
there exists $$q=[z_0,\dots,z_n]\in C$$ such that
$z_{\ell+2k+2}\neq 0$ and $z_{\ell+2k+3}=\dots=z_n=0$. Since 
$$\det(\partial_i\partial_jF(p))_{i,j=0,1}\neq 0,$$ we may assume that the same inequality holds for $q$. 
We may write 
$$
R_{\ell+k+1}=\alpha_{\ell+k+1} x_{\ell+2k+2}^3 + x_{\ell+2k+2}^2\cdot L_{\ell+k+1} + x_{\ell+2k+2} \cdot Q_{\ell+k+1} + R_{\ell+k+2}
$$
where $\alpha_{\ell+k+1}\in \mathbb C$, and $L_{\ell+k+1},Q_{\ell+k+1}, R_{\ell+k+2}\in \mathbb C[x_{I_{\ell,k+1}}]$ are homogeneous polynomials of degree $1$, $2$ and $3$ respectively. 



We first assume that $\alpha_{\ell+k+1}\neq 0$.  After replacing $x_{\ell+2k+2}$ by $x_{\ell+2k+2}- \frac 1 {3\alpha_{\ell+k+1}} L_{\ell+k+1}$, we may assume that $L_{\ell+k+1}=0$. Since $q\in V_F$, we get a contradiction by considering the minor $$(\partial_i\partial_j F(q))_{i,j=0,1,\ell+2k+2}.$$

We now assume that $\alpha_{\ell+k+1}= 0$. Since $z_{\ell+2k+2}\neq 0$ and $q\in V_F$ it follows that $L_{\ell+k+1}=0$ and that after a suitable change of coordinates, 
$Q_{\ell+k+1}\in \mathbb C[x_{\ell+1},x_{\ell+3},\dots,x_{\ell+2k+3}]$. 
We may write
$$Q_{\ell + k+1}=\beta_k x_{\ell+2k+3}^2 + x_{\ell+2k+3}\cdot \ell_{k} + M_k$$
where $\beta_k\in \mathbb C$ and  $\ell_{k},M_k\in \mathbb C[x_{\ell+1},x_{\ell+3},\dots,x_{\ell+2k+1}]$ are homogeneous polynomials of degree $1$ and $2$ respectively.
If $\beta_k\neq 0$ then, after a suitable change of coordinates, we may assume $\beta_k=1$ and $\ell_{k}=0$.  By considering the minor $$(\partial_i\partial_j F(q))_{i,j=0,\ell+2k+2,\ell+2k+3}$$ it follows that $C\subseteq \{x_{\ell+2k+3}=0\}$.
Thus, we may proceed to Step 4.

If $\beta_k=0$, then since $q\in V_F$ it follows that $\ell_k=0$. In case  $C$ is  contained in  $\{x_{\ell+2k+3}=\dots=x_n=0\}$ we are done, so we may assume that there exists a point 
$$q'=[z'_0,\dots,z'_n]\in C\cap \bigcap_{i\in J}\{x_i=0\}$$
such that $z'_0\neq 0$ and $z'_{\ell+2k+3}\neq 0$, where,  $J=I_{\ell,k+1}\setminus\{\ell+2k+3\}$.
 Proceeding as above, we may write
$$R_{\ell+k+2}=x_{\ell+2k+3}\cdot Q_{\ell+k+2} +R_{\ell+k+3},$$
where $Q_{\ell+k+2}\in \mathbb C[x_{\ell+1},x_{\ell+3},\dots,x_{\ell+2k+1},x_{\ell+2k+4}]$ and $R_{\ell+k+3}\in \mathbb C[x_{J}]$ 
are homogeneous polynomials of degree $2$ and $3$ respectively. We may write
$$
Q_{\ell + k+2}=\beta_{k+1} x_{\ell+2k+4}^2 + x_{\ell+2k+4}\cdot \ell_{k+1} + M_{k+1}
$$
where $\beta_{k+1}\in \mathbb C$ and  $\ell_{k+1},M_{k+1}\in \mathbb C[x_{\ell+1},x_{\ell+3},\dots,x_{\ell+2k+1}]$ are homogeneous polynomials of degree $1$ and $2$ respectively.  

If  $\beta_{k+1}= 0$ then $\ell_{k+1}=0$ because $q' \in V_F$. Denoting by $\mathcal H^i_F$ the $i$-th column of $\mathcal H_F$, it follows that  the vectors
$\mathcal H^{\ell+2}_F,\mathcal H^{\ell+4}_F\dots,\mathcal H^{\ell+2k+2}_F$ and $\mathcal H^{\ell+2k+3}_F$ are  linearly dependent. Thus, 
$\mathcal H_F$ does not have maximal rank which contradicts the assumptions.

Hence we have $\beta_{k+1}\neq 0$. After a suitable change of coordinates, we may assume that $\beta_{k+1}=1$ and $\ell_{k+1}=0$.  By considering the minor $$(\partial_i\partial_j F(q'))_{i,j=0,\ell+2k+3,\ell+2k+4}$$ it follows that $C\subseteq \{x_{\ell+2k+4}=0\}$. Thus we  first exchange $x_{\ell+2k+3}$ and $x_{\ell+2k+4}$, then we exchange $x_{\ell+2k+2}$ and $x_{\ell+2k+4}$. So we may  write
$$
R_{\ell+k+1}=x_{\ell+2k+2} \cdot (x_{\ell+2k+3}^2 + M_{k+1}) + R_{\ell+k+2}
$$
where $M_{k+1}\in \mathbb C[x_{\ell+1},x_{\ell+3},\dots,x_{\ell+2k+1}]$ is a quadric, $R_{\ell+k+2} \in \mathbb C[x_{I_{\ell,k+1}}]$ is a cubic form and $C \subseteq \{x_{\ell+2k+3}\}$. We also may write
$$
R_{\ell+k+2}= x_{\ell+2k+4}\cdot M_{k+2} +R_{\ell+k+3}
$$
where  $M_{k+2}\in \mathbb C[x_{\ell+1},x_{\ell+3},\dots,x_{\ell+2k+1}]$, $R_{\ell+k+3}\in \mathbb C[x_{I_{\ell,k+2}}]$  are  homogeneous polynomials of degree $2$ and 3 respectively. 

Moreover we have a point $$q'=[z'_0,\dots,z'_n]\in C\cap \bigcap_{i\in J}\{x_i=0\}$$
such that $z'_0\neq 0$ and $z'_{\ell+2k+2}\neq 0$, where $J=I_{\ell,k+1}\setminus\{\ell+2k+4\}$. Replacing $x_{\ell+2k+4}$ by $x_{\ell+2k+4}+\frac{z'_{\ell+2k+4}}{z'_{\ell+2k+2}}x_{\ell+2k+2}$ we get a point 
$$
q=[z_0,\ldots,z_n] \in C \cap \bigcap_{i \in I_{l,k+1}}\{x_i=0\}
$$ 

such that $z_0 \ne 0$, $z_{\ell+2k+2}\ne 0$ and we may proceed to Step 4.

\medskip 

\textbf{Step 4.}
We are assuming that
$$F=\sum_{i=0}^{\ell} G_i +\sum_{i=1}^k (x^2_{\ell+2i+1} + M_i)\cdot x_{\ell + 2i} + R_{\ell+k+1}$$
where $G_i,$ $M_i$ and $R_{\ell+k+1}$ satisfy \eqref{i_G}, \eqref{i_M} and \eqref{i_R} and $$C\subseteq \{x_{\ell+1}=x_{\ell+3}=\dots=x_{\ell+2k+1}=0\}.$$ By Step 3 we also have that
$$
R_{\ell+k+1}=x_{\ell+2k+2} \cdot (x_{\ell+2k+3}^2 + M_{k+1}) + R_{\ell+k+2}
$$
where $M_{k+1}\in \mathbb C[x_{\ell+1},x_{\ell+3},\dots,x_{\ell+2k+1}]$ is homogeneous of degree 2 and $C \subseteq \{x_{\ell+2k+3}=0  \}$. Moreover there is a point $q=[z_0,\ldots,z_n]$ such that $z_0 \ne 0$, $z_{\ell+2k+2}\ne 0$ and
$$
q \in C \cap \bigcap_{i \in I_{l,k+1}}\{x_i=0\}.
$$

We  show that we may assume
$$
M_{k+1}= \delta_{k+1} x_{\ell+1}^2
$$
where $\delta_k \in \mathbb C$.

Since $q \in C$ and $z_{\ell +2k+2} \ne 0$ we have $\det(\partial_i \partial_j F(q))_{i,j=0,1} = 0$. Considering the minors
$$
(\partial_i\partial_j F(q))_{i=0,h,\ell+2k+3}^{i=0,m,\ell+2k+3}
$$
for $h,m=1,\ldots,n$, $(h,m) \ne (\ell+2k+3,\ell+2k+3)$ we deduce that $\partial_h \partial_m F(q)=0$  and so, since by induction $M_{i}=\delta_i x_{\ell+1}$ for  $ i=1,\ldots k$, we have 
$$
M_{k+1}=\sum_{j=0}^{k} \gamma_k^j x^2_{\ell+2j+1},
$$
where $\gamma_k^j \in \mathbb C$. Since $M_{j}=\delta_j x_{\ell+1}$ for  $j=1,\ldots k$ to conclude it is enough  to replace $x_{\ell +2j}$ with $x_{\ell+2j}-\gamma_k^jx_{\ell+2k+2}$ for $j=1,\ldots,k$. In this way we get
$$
M_{k+1}= \delta_{k+1} x_{\ell+1}^2
$$
where $\delta_{k+1}=\gamma_k^0-\sum_{i=1}^k \gamma_k^i \delta_i$.

After replacing $k$ by $k+1$, we may repeat Step 3.  
\end{proof}

\begin{theorem}\label{t_plane}
Let $F\in \mathbb C[x_0,\dots,x_n]$ be a non-degenerate cubic form.

Then the set of points $p\in V_F$ such that $F(p)\neq 0$ is  a finite union of points, lines, plane conics and plane cubics. 
\end{theorem}
\begin{proof}
We may assume that there is an irreducible component $C\subseteq V_F$ such that $\dim C \ge1$ and  $F(p) \ne 0$ at the general point $p$ of $C$, otherwise we are done.  By Theorem \ref{t_special} we may write
$$
F=\sum_{i=0}^{\ell} G_i +\sum_{i=1}^k (x^2_{\ell+2i+1} + M_i)\cdot x_{\ell + 2i} + R_{\ell+k+1}
$$
where $G_i,$ $M_i$ and $R_{\ell+k+1}$ are as in Theorem \ref{t_special} and $$C\subseteq \{x_{\ell+1}=x_{\ell+3}=\dots=x_{\ell+2k+1}=0\}.$$

By the proof of Theorem \ref{t_special} we may also assume that for any $i=1,\ldots,k$ there is a point $q_i \in C$ such that $q_i \notin \{x_0=0\}$, $q_i \notin \{x_{\ell+2i}=0\}$ and $q_i \in \bigcap_{j=2i+1}^n \{x_{\ell+j}=0 \}$.  

We distinguish two cases: $C \subseteq \{x_1=0\}$ and $C \not\subseteq \{x_1 = 0\}$.

\medskip

If $C \subseteq \{x_1 = 0\}$ then  $\ell=0$. Let $z=[z_0,\ldots, z_n] \in C$ be a general point in $C$.

If $C \subseteq \{ x_{2k+2}=0\}$ then considering
$$
(\partial_i\partial_j F(z))_{i=0,1,2k+1}^{j=0,1,2k+1}
$$
we  immediately get a contradiction  because $\det(\partial_i\partial_jF(z))_{i,j=0,1}\neq 0$ and $z_{2k}\ne 0$.

So let $C \not\subseteq \{ x_{2k+2}=0\}$.  Then we may write
$$
R_{\ell+k+1}= M_{k+1}\cdot x_{\ell+2k+2}+R_{l+k+2}
$$ 
as in \eqref{i_k} of Theorem \ref{t_special}.
Assume that $k >2$. Then we have 
\begin{align*}
&\det(\partial_i\partial_jF)_{i=0,1,2k+1}^{j=0,3,2k+1}=  \\  &=6x_0\cdot(2\gamma_{1,3}x_{2k}x_{2k+2}+\gamma_{1,3}\gamma_{2k+1,2k+1}x^2_{2k+2}-\gamma_{1,2k+1}\gamma_{3,2k+1}x^2_{2k+2}+ Q)
\end{align*}
where $Q \in \C[x_1,\ldots,x_n]$ is a quadratic form such that $C\subseteq \{Q=0\}$ (because  $C \subseteq\bigcap_{i\in I_{\ell,k+1}}\{x_i=0 \}$) and where $\gamma_{i,j}$ is the coefficient of $x_{2k+2}$ in $\partial_i\partial_j F$. Note that $\gamma_{1,3} \ne 0$ (because $\partial_3\partial_3 F(z) \ne 0$, being this last inequality true for $q_2$).

Since $z_0 \ne 0$ and $z_{\ell+2k} \ne 0$ we conclude that 
$$C \subset\{ 2\gamma_{1,3}x_{2k}+(\gamma_{1,3}\gamma_{2k+1,2k+1}-\gamma_{1,2k+1}\gamma_{3,2k+1})x_{2k+2}=0 \},$$ which contradicts the fact that $q_{k} \in C$.  Hence we conclude that $k \le 2$. Now it is easy to see that $C$ is a line or a plane conic.

\medskip

Assume now that $C \not\subseteq \{x_1=0\}$. Then $\ell\ge 1$. Note that for $j=3,\ldots,n$ we have $\partial_1\partial_jF=0$, hence for a general point $z=[z_0,\ldots,z_n] \in C$, for $h=2,\ldots,n$ and for $m=3,\ldots,n$ we may consider 
$$
(\partial_i\partial_j F(z))_{i=0,1,h}^{j=0,1,m}
$$
to conclude that $\partial_h\partial_mF(z)=0$ (because $\det(\partial_i\partial_jF(z))_{i,j=0,1}\neq 0$). This implies easily that we may assume $k=0$. By Step 2 of the proof of Theorem  \ref{t_special}  for any $i=1,\ldots,\ell$ there is a point $p_i \in C$ such that $p_i \notin \{x_0=0\}$, $p_i \notin \{x_{i}=0\}$ and $p_i \in \bigcap_{j=i+1}^n \{x_j=0 \}$.

\medskip
Assume first that $C \subseteq \{x_{\ell+2}=0 \}$ so we may write 
$$
F=\sum_{i=0}^{\ell} G_i + R_{\ell+1}
$$
where $G_i \in \C[x_i,x_{i+1}]$, $R_{\ell+1} \in \C[x_{\ell+1},\ldots,x_n]$ are cubic forms and $C \subseteq \bigcap_{i=\ell+1}^n\{x_i=0\}$.

Suppose  that $\ell>2$. Since $\partial_{3}\partial_{3}F(p_2)=0$, $\partial_2\partial_3F(p_2)=0$ and $\partial_{3}\partial_{3}F(p_3)=0$ we see that the monomials $x_2x_3^2$, $x_2^2x_3$ and $x_3^3$ do not appear in $F$.  
The same holds for $x_3x_4^2$ and $x_3^2x_4$ which gives a contradiction. Hence $\ell \le 2$ and it is easy to conclude.

\medskip

If $C\not\subseteq \{x_{\ell+2}=0 \}$ then we may write
$$
F=\sum_{i=0}^{\ell} G_i +x^2_{\ell+1} \cdot x_{\ell + 2} + R_{\ell+1}.
$$
where $G_i \in \C[x_i,x_{i+1}]$ and $R_{\ell+1} \in \C[x_{I_{\ell,1}}]$.

Suppose $\ell \ge 2$. Since $\partial_{\ell+1}	\partial_{\ell+1}F(p_\ell)=0$ we see that $x_{\ell+1}^2x_\ell$ does not appear in $F$ and this implies, considering $\partial_{\ell+1}\partial_{\ell+1}F(z)$, that also $x_{\ell+1}^2x_{\ell+2}$ does not appear in $F$, which is a contradiction. Thus $\ell <2$ and we are done.
\end{proof}

\begin{remark}
Note that in general $V_F$ might contain surfaces, e.g. if
$$F(x_0,\dots,x_n)=x_0^3+x_0x_{1}^2 + x_{1}\cdot\sum_{i=2}^{n} x_i^2$$
then $\dim V_F=n-2$. 
\end{remark}

\medskip

Our goal is now to improve Theorems \ref{t_special} and \ref{t_plane}   and characterise those cubic forms $F$ such that $V_F$ contains a curve $C$ such that $C\nsubseteq \{F=0\}$.  To this end, we restrict to the case of cubic forms  
with non-zero discriminant. 

\begin{corollary}\label{c_Csub}
Let $F\in \mathbb C[x_0,\dots,x_n]$ be a non-degenerate cubic form such that
$$
F=ax_0^3 + bx_0^2x_1+G(x_1,\ldots,x_n).
$$
  Let $C \subseteq V_F$ be positive dimensional irreducible variety such that $p=[1,0,\ldots,0] \in C$ and assume that at least one of the following properties holds:
\begin{enumerate}
\item $C \subseteq \{x_1=0\}$;
\item $C\subseteq \{F=0\}$.
\end{enumerate}
Then $\Delta_F=0$.
\end{corollary}
\begin{proof}
We first assume that $C\subseteq \{x_1=0\}$.
By the proof of Theorem \ref{t_plane}, we may write 
$$
F=x_0^3 + x_0x_1^2 +  (x_3^2 + \delta_1x_1^2)x_2 + R(x_1,x_3,x_4,\ldots,x_n),
$$
for some $\delta_i \in \mathbb C$ and $R\in \mathbb C[x_1,x_3,x_4,\dots,x_n]_3$. It follows that the hypersurface $\{F=0\}\subseteq \mathbb P^n$ is singular at the point $[0,0,1,0,\dots,0]$ and in particular $\Delta_F=0$, as claimed.

\medskip

We now suppose that $C\subseteq \{F=0\}$ and $C \not\subseteq \{ x_1=0\}$.
Since $[1,0,\dots,0]\in C$ we may write
$$
F= bx_0^2x_1+ c_1 x_1^3 + Lx_1^2 + Qx_1 +R
$$
where  $b,c_1\in \C$ and $L,Q,R \in \mathbb C[x_2,\ldots,x_n]$ are homogeneous polynomials of degree 1,2 and 3 respectively. Since $F$ is non-degenerate, we have that $b\neq 0$.

%
%

\medskip

After a change of coordinates in $(x_1,x_2,\ldots,x_n)$ we may assume that there exists a point $q=[q_0,q_1,0,\ldots,0] \in C$ such that $q_0,q_1 \ne 0$ and that $L=c_2x_2$ for some $c_2\in \mathbb C$. 
Note that since $C\subseteq \{F=0\}$, it follows that $C$ is not  a line. 
Furthermore, since $q\in V_F$ we may assume that $Q=c_3x_2^2$ for some $c_3\in \mathbb C$ and we may write
$$
F= bx_0^2x_1+ c_1 x_1^3 + c_2x_1^2x_2 + c_3x_1x_2^2 + c_4x_2^3 +R_1
$$
where $c_4\in \mathbb C$ and  $R_1 \in \mathbb C[x_2,\ldots,x_n]_3$ is such that  the monomial $x_2^3$ does not appear in $R_1$. It is easy to see that $\partial_i\partial_j F(z)=0$ for $i=2,\ldots,n$, $j=2,\ldots,n$, with $(i,j) \ne (2,2)$ and  $z \in C$.    If $ C \subseteq \{ x_2=0\}$ then, after a change of coordinates in $(x_3,\ldots,x_n)$, we may assume that there is a point $r=[r_0,r_1,0,r_3, 0, \ldots, 0] \in C$	such that $r_3 \ne 0$.
It follows that
$$
R_1= \alpha x_2^2x_3 + R_2(x_2,x_4,\ldots,x_n),
$$
for some $\alpha \in \mathbb C$ and $R_2 \in \C[x_4,\ldots,x_n]_3$. In particular, $[0,0,0,1,0\ldots,0]$ is a singular point of $\{F=0\}\subseteq \mathbb P^n$. Thus, $\Delta_F=0$, as claimed. 

Thus, we may assume that  $C \not\subseteq \{x_2=0 \}$ and  that there is a point $s=[s_0, s_1, s_2, 0, \ldots,0]$ such that $s_2 \ne 0$. Since $\partial_i\partial_j F(s)=0$ for $i=2,\ldots,n$, $j=2,\ldots,n$, with $(i,j) \ne (2,2)$, it follows that $R_1$ does not depend on $x_2$. Thus, $\partial_i\partial_j F(z)=0$  for any $i,j\ge 3$ and $z \in C$. Lemma \ref{l_trivialh} implies that  $C$ is contained in the plane $\Pi=\{ x_3= \ldots=x_n=0\}$. 
Let $F_1$ be the restriction of $F$ to $\Pi$. Since $C\subseteq \{F=0\}$, it follows that if $[x_0,x_1,x_2,0,\dots,0]\in C$ then $F_1(x_0,x_1,x_2)=0$ and $\mathcal H_{F_1}(x_0,x_1,x_2)=0$. Thus $C$ is a line, which gives a  contradiction. 
\end{proof}

\begin{corollary}\label{c_Cline}
Let 
$$
F(x_0, \ldots, x_n)=ax_0^3 + x_0^2(bx_1+cx_2) + G(x_1,\ldots,x_n)
$$
be a non-degenerate cubic form with integral coefficients such that  $b \ne 0$.
Assume that the line  $C=\{x_2=x_3=\ldots=x_n=0\}$ is contained inside  $V_F$.  

Then there exists $T=(t_{ij})_{i,j=0,\ldots,n} \in \SL(n+1,\mathbb Q)$ such that
$$
T\cdot F =ax_0^3 + bx_0^2x_1 + c_1 x_1^3 + R(x_2,\ldots,x_n)
$$	
where  $c_1 \in \ZZ$ and $R \in \mathbb Q[x_2,\ldots,x_n]$ is a cubic form. Moreover we may choose $T$ such that  $t_{00}=t_{11}=1$, $t_{0i}=t_{i0}=0$ for $i=1,\ldots,n$,  $t_{ij}=0$ for $i=2,\ldots,n$ and $j=1$  
\end{corollary}

\begin{proof} 
After replacing $x_1$ by $x_1 - cx_2/b$, we may write 
$$F=ax_0^3 + bx_0^2x_1+ c_1 x_1^3 + Lx_1^2 + Qx_1 +R$$
where $c_1 \in \ZZ$ and $L,Q,R \in \mathbb Q[x_2,\ldots,x_n]$ are homogeneous polynomials of degree 1,2 and 3 respectively. After a change of coordinates in $(x_2,\ldots,x_n)$ we may also assume that $L= c_2 x_2$, for some $c_2\in \mathbb Q$. 
Let $q=[0,1,0\ldots,0] \in C$. We distinguish two cases: $c_1 \ne 0$ and $c_1 =0$.

\medskip

If $c_1 \ne 0$ then,  since $b\neq 0$ and $\rk \mathcal H_F(q)\le 2$, we see that $Q= c_3 x_2^2$ for some $c_3\in \mathbb Q$ and $$|(\partial_i\partial_j F(q))_{i=1,2}|=0.$$ 
It follows that $|(\partial_i\partial_j F(z))_{i=1,2}|=0$ for any $z \in C$. 
Since  
$$|(\partial_i\partial_j F(z))_{i,j=0,1,2}|=0,$$ we have that $c_2=c_3=0.$   Thus, $L=Q=0$ and the claim follows. 

\medskip

If $c_1 =0$ then since  $b\neq 0$ and $\rk \mathcal H_F(q)\le 2 $, it follows that $c_2=0$.
 Since $\rk \mathcal H_F(z)\le 2$ for any $z \in C$, we have $Q=0$ and, again, the claim follows. Note that in this case, we have $\Delta_F =0$.
\end{proof}

\begin{corollary}\label{c_Ccubic}
Let 
$$
F(x_0, \ldots, x_n)=ax_0^3 + x_0^2(bx_1+cx_3) + G(x_1,\ldots,x_n)
$$
be a non-degenerate cubic form with integral coefficients with $b,c\in \mathbb Z$ and $G\in \mathbb Z[x_1,\dots,x_n]$ such that $b\neq 0$ and  $\Delta_F\neq 0$. 
Let $C \subseteq V_F$ be a positive dimensional  irreducible variety  such that $C \not\subseteq\{F=0\}$ and $p=[1,0,\ldots,0] \in C$. Assume that $C$ contains infinitely many rational points. Assume moreover that $C \subseteq\Pi=\{x_3=\ldots=x_n=0\}$ and $C$ is not a line. 

Then there exists $T=(t_{ij})_{i,j=0,\ldots,n} \in \SL(n+1,\mathbb Q)$,  $R \in \ZZ[x_1,x_2]_3$ and $S \in \mathbb Q[x_3,\ldots,x_n]_3$  such that:
\begin{enumerate}
\item  $t_{00}=1$, $t_{i0}=t_{0i}=0$ for $i=1,\ldots,n$,  $t_{ij}=0$ for $i=3,\ldots,n$ and $j=1,2$, $(t_{ij})_{i,j=0,1,2} \in \SL(3,\ZZ)$ and 
\item $T\cdot F =ax_0^3 + bx_0^2x_1 + R(x_1,x_2) + S(x_3,\ldots,x_n).$
  

\end{enumerate}	
\end{corollary}

\begin{proof}
We may assume that there is a point $q=[z_0,1,0,\ldots,0] \in C$ such that $z_0 \ne 0$. Indeed, since $C$ is not a line, there exists  $m \in \ZZ$ such that $\{mx_1+x_2=0\}\cap \Pi$ intersect $C$ in a point $[z_0,1,-m,0,\ldots,0]$ with $z_0 \ne 0$. After replacing $x_2$ with $x_2 + mx_1$, we may assume that $m=0$. 

In addition, after replacing  $x_1$ with $x_1- c/bx_3$, we may assume that $c=0$. Thus, we may write
$$
F=ax_0^3 + bx_0^2x_1 + c_1x_1^3 + c_2x_1^2x_2+c_3x_1x_2^2 + c_4x_2^3+x_1^2 L + x_1 Q + S
$$
where $c_i \in \ZZ$ and $L\in \mathbb Q[x_3,\dots,x_n]$, and $Q,S \in \mathbb Q[x_2,\ldots,x_n]$ are homogeneous polynomials of degree 1,2 and 3 respectively such that the coefficient of $x_2^2$ in $Q$ and the coefficient of $x_2^3$ in $S$ are zero.

\medskip

If $c_2 \ne 0$ then, after  replacing $x_2$ with $x_2 -L/c_2$, we may assume $L=0$. Since $b\neq 0$ and 
$q\in V_F$, it follows that $Q=0$.  Now considering a general point $z \in C\subseteq \{x_3=\ldots=x_n=0\}$, we see that $\partial_i\partial_j S(1,0,\dots,0)=0$ for all $i,j\ge 2$. As in the proof of Lemma \ref{l_trivialh}, it follows that $S$ does not depend on $x_2$. Thus, (2) holds.

\medskip
Assume now that $c_2=0$ and $L=0$. Then the Hessian  of the quadric $c_3x_2^2+Q$ has rank  not greater than $1$, which means that
$$
c_3x_2^2+Q=c_3(x_2+L_1)^2
$$
for some $L_1 \in \mathbb Q[x_3,\cdots,x_n]$ of degree $1$. Hence, replacing $x_2$ with $x_2 - L_1$  we may assume that $Q=0$. As in the previous case, it follows that $S$ does depend on $x_2$. Thus, (2) holds.

\medskip

Finally assume  that $c_2=0$ and $L \ne 0$. Acting on $(x_3,\ldots,x_n)$ with $\SL(n-2,\mathbb Q)$ we may 
write $L=\alpha x_3$, where $\alpha \ne 0$. 
In particular,  $\partial_3\partial_1F(q)\neq 0$. It follows that the first two columns $\mathcal H^0_F(q)$ and $\mathcal 
H^1_F(q)$ of $\mathcal H_F(q)$ are linearly independent, which implies that $c_3=0$. Considering now a 
general point in $ C\subseteq \{x_3=\dots=x_n=0\}$, we  see that $c_4=0$. and that the only monomial which appears in $x_1Q+S$ with non-zero coefficient and which contains $x_2$ is  
 $x_2x_3^2$. Since $[0,0,1,0,\ldots,0]$ is a singular point of the hypersurface $\{F=0\}
\subseteq \mathbb P^n$, it follows that 
$\Delta_F=0$, a contradiction.
\end{proof}

\subsection{Binary and ternary cubic forms}
We now study the possible reduced forms of a non-degenerate binary or ternary cubic form. We show that if $F$ is a binary cubic form, it admits only  finitely many non-equivalent reduced forms (cf. Proposition \ref{p_binary}). On the other hand, if $F$ is a ternary cubic form, then the same result holds with the extra assumption that the discriminant $\Delta_F$ is non-zero (cf. Proposition \ref{p_b2=3}).  Example \ref{ex_D=0} shows that this assumption is necessary.

We first recall the following known result: 

\begin{proposition}\label{p_finite}
 Let $\Delta \ne 0$ be an integer. Then there exist 
$$F_1,\dots,F_k\in \mathbb Z[x_0,x_1,x_2]_3\qquad \text{(resp. $\mathbb Z[x_0,x_1]_3$)}$$ such that if $F\in \mathbb Z[x_0,x_1,x_2]_3$  (resp. $\mathbb Z[x_0,x_1]_3$) is such that $\Delta_F=\Delta$, then there exists $i=1,\dots,k$ and $T\in \SL(3,\mathbb Z)$ (resp. $\SL(2,\mathbb Z)$) such that $F=T\cdot F_i$. 
\end{proposition}
\begin{proof}
See \cite[Proposition 7]{OV95}.
\end{proof}
\medskip


\begin{lemma}\label{l_b2=2}
Let $$F(x,y)= ax^3+bx^2y+cy^3\in \mathbb Z[x,y]$$ be a binary cubic form with integral coefficients and  such that $c \ne 0$. 

Then there are  finitely many pairs
$$
(a_i,b_i)\in \mathbb Z^2 \qquad i=1,\dots,k
$$  
such that if $(a',b',cy^3)$ is a reduced triple associated to $F$ (cf. Definition \ref{d_reduced}) then  $a'=a_i$ and $b'=b_i$ for some $i \in \{1,\dots,k\}$.  
\end{lemma}

\begin{proof}
%
%
Assume that 
$T=(t_{i,j})_{ij=0,1}\in \SL(2,\mathbb Z)$ is such that 
$T\cdot F$ is in reduced form $(a',b',cy^3)$, for some $a',b'\in \mathbb Z$. 

Note that  $F(t_{01},t_{11})=c$ and, since $c \ne 0$, the equation $F(x,y)=c$ defines a smooth affine plane curve of genus 1. Thus, by Siegel's Theorem \ref{t_siegel}, it only admits finitely many solutions. Thus, we may assume that $t_{01}$ and $t_{11}$ are fixed. 
Since $\det T=1$ and since the coefficient of $xy^2$ is  zero, 
we get the linear system  in $t_{00}$ and $t_{10}$:
$$\begin{cases}
 1 =&t_{11} t_{00}- t_{01}t_{10} \\
0= &(3a t_{01}^2\ +2b t_{01}t_{11})t_{00} + (b t_{01}^2 + 3c t_{11}^2)t_{10}. 
\end{cases}$$
Note that the determinant of the system is equal to $3F(t_{01},t_{11})=3c\neq 0$. Thus, 
the system admits exactly  one solution and the claim  follows. 
\end{proof}

%
%
%



\begin{proposition}\label{p_binary}
Let $$F(x,y)= ax^3+bx^2y+cy^3 \in  \Z[x,y]$$ be a binary integral cubic form with $c \ne 0$. 

Then there are  finitely many triples 
$$
(a_i,b_i, c_i)\in \Z^3\qquad i=1,\dots,k
$$  
such that $c_i \ne 0$ and if $(a',b',c'y^3)$ is a reduced triple associated to $F$ (cf. Definition \ref{d_reduced}) then $a'=a_i$, $b'=b_i$ and $c'=c_i$ for some $i \in \{1,\dots,k\}$.
\end{proposition}

\begin{proof}
By Lemma \ref{l_b2=2}, it is enough to show that there are only finitely many  $c_1,\dots,c_k \in \Z$ such that if $T\in \SL(3,\mathbb Z)$ is such that $T\cdot F$ is in reduced form $(a',b',c'y^3)$, with $c'\neq 0$, then $c'=c_i$ for some $i \in \{1,\dots,k\}$. 

If the discriminant  $\Delta_F= 4b^3c +27a^2c^2$ of $F$ is not zero, then $c' | \Delta_F$ and the claim follows.

Thus, we may assume that $\Delta_F=0$. We may also assume that $a,b$ and $c$ do not have a common factor, otherwise we just consider the cubic form obtained by dividing by the  common factor. Suppose that $T=(t_{ij})_{i,j=0,1}$. 
Then,
\begin{align}
\label{F_a} a&= a't_{00}^3 + b't_{00}^2t_{10}+c't_{10}^3, \\
\label{F_b} b&=3a't_{00}^2t_{01} + b't_{00}^2t_{11} + 2b't_{00}t_{01}t_{10} + 3c't_{10}^2t_{11}, \\
\label{F_0} 0&=3a't_{00}t_{01}^2 + b't_{01}^2t_{10} + 2b't_{00}t_{01}t_{11} + 3c't_{10}t_{11}^2, \\
\label{F_c} c&= a't_{01}^3 + b't_{01}^2t_{11}+c't_{11}^3
\end{align} 
and $GCD(a',b',c')=1$.

Let $p$ be a prime factor of $c'$ such that $p \ne 2,3$ and let $\alpha$ be a positive integer such that $p^\alpha | c'$.  Then, since $\Delta_F=0$, it follows that   $p^{\rup\alpha/3.}$ divides $b'$. 
By \eqref{F_0}, and since  $\gcd(t_{00},t_{01})=1$, we have that either $p^{\rup \alpha/3.}$ divides $t_{00}$ or $p^{\rup \alpha/6.}$ divides $ t_{01}$. In the first case, \eqref{F_a} implies that  $p^{\alpha}$ divides  $a$, and in the second case,  \eqref{F_c} implies that  $p^{\rup \alpha/2.}$ divides $c$ . Since $a,c\neq 0$ are fixed, it follows that  $p^{\alpha}$ is bounded. A similar argument holds for the powers of $2$ and $3$. Hence $c'$ is bounded, as claimed.
\end{proof}

We now consider ternary cubic forms:

\begin{proposition}\label{p_Dne0}
Let $R$ be a ring which is finitely generated over $\mathbb Z$ and let  $F \in R[x,y,z]$ be a cubic form with non-zero discriminant $\Delta_F$.  Let $G(y,z)=dy^3+z^3$ for some non-zero $d\in R$ and assume that $F$ is in reduced form $(a,(b,c),G)$ for some pair $(a,(b,c))\in R\times R^2$. 

Then there are finitely many  pairs
$$
(a_i,(b_i,c_i)) \in R \times R^2\qquad i=1,\dots,k
$$ 
  such that  if $(a',(b',c'),G)$ is a reduced triple associated to $F$ (cf. Definition \ref{d_reduced}) then  $a'=a_i, b'=b_i$ and $c'=c_i$ for some $i \in \{1,\dots,k\}$.  
\end{proposition}
\begin{proof}
Assume that $T\in \SL(3,R)$ is  such that $T\cdot F$ is in reduced form $(a',(b',c'), G)$. 
The invariants $S_F$ and $T_F$ (cf. Subsection \ref{s_cf} and \cite[ 4.4.7 and 4.5.3]{Sturmfels93}) have the form 
$$
S_F=dbc\quad \mbox{ and } \quad T_F=27a^2d^2+4b^3d + 4c^3d^2. 
$$

We first assume that $S_F \ne 0$ and we consider the curve $C \subseteq\mathbb P^3$ given by the ideal 
$$
I=(S_F x_3^2-dx_1x_2, T_F x_3^3-27d^2x_0^2x_3-4dx_1^3-4d^2x_2^3).
$$
 
We claim that the  points $[a',b',c',1]\in C$, with $a',b',c'\in R$ are in finite number and hence the claim follows. 

Note that the first equation define a cone over a conic with vertex the point $q=[1,0,0,0]\in C$. If we blow-up the  point $q$, then it is easy to check  the strict transform $\tilde C$ of the curve $C$ is a connected smooth curve of genus $3$.
Thus, the claim follows by Siegel's Theorem \ref{t_siegel}.

\medskip

We now assume that $S_F=0$. Then, $b'=0$ or $c'=0$. 
 Assume that $c'=0$. Then the pair $(a',b')$ corresponds to an $R$-integral point in  the affine plane curve, defined by the equation 
$$
27x_0^2d^2+4x_1^3d-T_F=0.
$$
Since, by assumption $\Delta_F\neq 0$, we have that  $T_F \ne 0$. Thus, Siegel's Theorem \ref{t_siegel} implies the claim. The case $b'=0$ is similar.  
\end{proof}

\begin{remark}
Note that if $F\in R[x,y,z]$ is a cubic form such that $\Delta_F=0$ and $S_F=0$, and $C$ is the curve defined in the proof of Proposition \ref{p_Dne0}, then $C$ is a rational curve.
\end{remark}

As a consequence of the previous result we  obtain the following:

\begin{proposition}\label{p_b2=3}
Let $F\in \ZZ[x,y,z]$ be a cubic form with non-zero discriminant  $\Delta_F$. 

Then there are  finitely many triples 
$$(a_i,B_i,G_i)\in \mathbb Z\times \mathbb Z^2\times \mathbb Z[y,z]_3\qquad i=1,\dots,k$$ 
such that any reduced triple associated to $F$ is 
equivalent to $(a_i,B_i,G_i)$ over $\mathbb Z$, for some $i \in \{1,\dots,k\}$ (cf. Definition \ref{d_reduced}).
\end{proposition}

\begin{proof}
Let $T\in \SL(3,Z)$ such that $T\cdot F$ is in reduced form $(a,B,G)$ for some $a\in \mathbb Z$, $B\in \mathbb Z^2$ and $G\in \mathbb Z[y,z]$ cubic form. 
Lemma \ref{l_disc} implies  that $\Delta_{G}$ divides $\Delta_F$. Thus, $\Delta_G\neq 0$ and we may assume that its  value  is fixed, and, by  
 Proposition \ref{p_finite}, we may assume that  $G$ is also fixed, up to the action of $\SL(2,\mathbb Z)$.
 
Let $d=\sqrt{\frac{\Delta_F}{27}} $.
After possibly replacing the ring of integers $\mathbb Z$ by a finitely generated ring $R$ over $\mathbb Z$, we may assume, up to a $\SL(2,R)$-action, that 
 $$G(y,z)=dy^3+z^3.$$
 
 Thus, the claim follows from  Proposition \ref{p_Dne0}.
  \end{proof}

\medskip

Note that Proposition \ref{p_Dne0} does not hold if the discriminant of $F$ is zero, as the following example shows:

\begin{example}\label{ex_D=0}
Let 
$$
F=ax^3 + bx^2y+ x^2z -3y^2z
$$
where $a,b \in \ZZ$. Note that $\Delta_F=0$, since $[0,0,1]$ is a singular point for $\{F=0\}$. Consider the Pell's equation
\begin{align}
\label{Pe1} s^2-3t^2=1. 
\end{align}
For any solution $(\alpha, \beta) \in \ZZ^2$ of \eqref{Pe1}, we define the matrix 
$$M=\begin{pmatrix}
\alpha & 3\beta& 0\\
\beta &\alpha &0 \\
m_{31} & m_{32} & 1
\end{pmatrix}$$
where $m_{31}=\beta(3b\beta^2+9a\alpha\beta+2b\alpha^2)$ and $m_{32}=3\beta^2(3a\beta+b\alpha)$.

Then $M \in \SL(3, \ZZ)$ and 
$$
M\cdot F(X,Y,X)= AX^3 + BX^2Y + X^2Z -3Y^2Z.
$$
where
$$
A=3b\alpha^2\beta+3b\beta^3+a\alpha^3+9a\alpha\beta^2 \mbox{ and } B=9a\beta^3+9b\alpha\beta^2+9a\alpha^2\beta+b\alpha^3.
$$

Since \eqref{Pe1} has infinitely many  integral solutions,  it follows that there are infinitely many ways to write $F$ in reduced form.  
\end{example}

In the example above, $\{F=0\}$ defines an irreducible cubic with a node.  Note that such cubics can be realised  as the cubic form associated to a smooth  threefold (the existence of such a threefold was asked in   \cite[Proposition 21]{OV95}):

\begin{example}
Let $W=\mathbb P^3$, $h$ the hyperplane class and $C$ a line. Note that $\deg N_{C / W}=2$. Let $\pi: X \to W$ be the blow-up of $W$ along $C$ and define $H=\pi^*h$. Let $\{L_1,L_2\}$ be the basis of $H^2(X,\ZZ)$ given by
$$
L_1=H \quad \mbox{ and } L_2=H-E
$$
where $E$ is the exceptional divisor  of $\pi$. 
The intersection cubic form on $H^2(X,\ZZ)$ is
$$
G(y,z)=(yL_1+zL_2)^3=y^3+3y^2z.
$$

Let  $C' \subseteq\mathbb P^3$ be a line which meets  $C$ transversally
in one point 
and let $D$ be the strict transform of $C'$ in $X$. Then  $D\equiv H^2 -H\cdot E$ and blowing-up $X$ along $D$  we get a threefold $Y$ with associated cubic form
$$
F(x,y,z)=x^3 -3(y+z)x^2+y^3+3y^2z.
$$
Note that $\{F=0\}\subseteq \mathbb P^2$ defines an irreducible  cubic with a node and in particular $\Delta_F=0$. 
\end{example}

\medskip 

\subsection{General cubic forms}\label{s_general}

We now combine the previous results to give a proof of Theorem \ref{t_b2=n}. We begin with the following:

\begin{lemma}\label{l_fp}
Let $F\in \Z[x_0,\dots,x_n]$  be a non-degenerate cubic form and let $p\in V_F$ such that $F(p)\neq 0$. 

Then there are  finitely many triples
$$
(a_i,B_i,G_i)\in \mathbb Z\times \mathbb Z^n \times \mathbb Z[x_1,\dots,x_n]_3\qquad i=1,\dots,k
$$
such that for all $T\in \SL(n+1,\mathbb Z)$  such that $T \cdot p = [1,0,\dots,0]$ and $T\cdot F$ is in reduced form, we have that $T\cdot F$  is equivalent to $(a_i,B_i,G_i)$ over $\mathbb Z$ for some $i \in \{1,\dots, k\}$ (cf. Definition \ref{d_reduced}). 
\end{lemma}

\begin{proof}
We may assume that $p=[1,0,\dots,0]$ and that $F=(a,b,G)$ is in reduced form, for some $a\in \mathbb Z$, 
$B\in \mathbb Z^n$ and $G\in \mathbb Z[x_1,\dots,x_n]_3$. We consider all the matrices $T\in \SL(n+1,\mathbb Z)$ such that
$T\cdot p=p$ and $T\cdot F=(a_T,b_T,G_T)$ is in reduced form, for some $a_T\in \mathbb Z$, 
$B_T\in \mathbb Z^n$ and $G_T\in \mathbb Z[x_1,\dots,x_n]$.

If we write $ T= (t_{ij})_{i,j=0,\ldots,n}$ with $t_{ij}\in \mathbb Z$,
then, since $T\cdot p=p$, we have $t_{i0}=0$ for $1\le i \le n$. Thus,  $t_{00}=\pm 1$ and in particular $a_T=\pm a$. 

By considering the action of $\SL(n,\ZZ)$ over $(x_1,\ldots,x_n)$, we may assume that $B=(b_1,0,\ldots,0)$ and that, for each $T$, $B_T=(b^T_1,0\ldots,0)$, with $b_1,b_1^T\in \mathbb Z$. 
Note that, by the assumption on $F$, we have that  $a$ and $b_1$ cannot be both zero.

By looking at the coefficients of $x_0^2x_i$ and $x_0x_i^2$, we obtain the equations
\begin{equation}\label{eq_a}
\begin{aligned} 
&3at_{0i}+b_1t_{1i}=0 \qquad\text{for $i=2,\dots,n$ and}\\ 
& 3at_{0i}^2 + 2b_1t_{0i}t_{1i}=0 \qquad\text{for $i=1,\dots,n$} . 
\end{aligned}
\end{equation}
 We now consider three cases.

\medskip

If $b_1=0$ then $a\neq 0$ and  \eqref{eq_a} implies that $t_{0i}=0$ for $i=1,\ldots,n$. In particular, $T\cdot F$ is equivalent to $F$.

\medskip

If $a=0$ then $b_1\neq 0$ and \eqref{eq_a} implies that $t_{1i}=0$ for $i=2,\ldots,n$. In particular, $t_{11}=\pm 1$. By looking at the coefficients of $x_0x_1x_i$ for $i=1,\ldots,n$, we get the equations
$$
b_1t_{0i}t_{11}=0.
$$
Thus $t_{0i}=0$ for $i=1,\ldots,n$ and, as in the previous case, we obtain that $T\cdot F$ is equivalent to $F$. 

\medskip 

Finally if $a,b \ne 0$ then \eqref{eq_a} implies that $t_{0i}=t_{1i}=0$ for $i=2,\dots,n$. In particular, $t_{11}= \pm 1$. By  \eqref{eq_a}, it follows that  $t_{01}$ can only acquire finitely many values. Thus, under these assumptions on $T$, it follows that there are only finitely many non-equivalent reduced form $T\cdot F$ over $\mathbb Z$, as claimed. 
\end{proof}

In the next Lemma we show that under the action of the transformations given by Corollaries \ref{c_Cline} and \ref{c_Ccubic} we may control the last part of a reduced form.

\begin{lemma}\label{l_G}
Let $s \in \{1,2\}$ and let $F,F_1 \in \mathbb Q[x_0,\ldots,x_n]$ be non-
degenerate cubic forms such that
$$
F=ax_0^3+ bx_0^2x_1+R(x_1,x_s) + H(x_{s+1},\ldots,x_n)
$$ 
and
$$
F_1=a_1x_0^3+ b_1x_0^2x_1+R_1(x_1,x_s) + H_1(x_{s+1},\ldots,x_n)
$$ 
where $b,b_1 \ne 0$ and $R,R_1, H, H_1$ are cubic forms.

Assume that there exists $T =( t_{hk})_{h,k=0,\ldots,n}\in \SL(n+1,\mathbb Q)$   such that $T \cdot F= F_1$, $t_{hk}=0$ for $h=s+1,\ldots,n$ and $k= 0,\ldots,s$ and $\det (t_{hk})_{h,k=0,\ldots,s}=1$, i.e.\
$$
T=\begin{pmatrix}
S & * \\
0 & *  
\end{pmatrix}
$$
with $\det S=1$.

Then there exists $P\in \SL(n-s,\mathbb Q)$ such that $P\cdot H=H_1$.
\end{lemma}

\begin{proof}
We prove the case $s=2$, the case $s=1$ is similar and easier.

\medskip
We will show that $t_{hk}=0$ for $h=0,1,2$ and $k=3,\ldots,n$, which implies the claim.

Let $S=(t_{hk})_{h,k=0,1,2}$ and define $\overline{T}=(\overline{t_{hk}})_{h,k=0,\ldots,n} \in \SL(n+1,\mathbb Q)$ as 
$$
\overline{T}=\begin{pmatrix}
S^{-1} & 0 \\
0 & I_{n-2}  
\end{pmatrix}
$$
where $I_{n-2} \in \SL(n-2,\mathbb Q)$ is the identity matrix. 

If $M=(m_{ij})_{i,j=0,\ldots,n}=\overline{T} \cdot T$ and  $\overline{F_1}=M\cdot F$, then $\overline{F_1}$ is in reduced form with associated triple $(a,(b,0),R+H_1)$.  In addition
$$\begin{aligned}
(m_{hk})_{h,k=0,1,2}&= I_3, \quad\text{and} \\
   (m_{hk})_{h=3,\ldots,n}^{k=0,1,2}&=0. 
   \end{aligned}
$$

\medskip

We want to show that $m_{hk}=0$ for $h=0,1,2$ and $k=3,\ldots,n$. Since $S$ is invertible, it follows that  $t_{hk}=0$ for $h=0,1,2$ and $k=3,\ldots,n$, as claimed.

We assume first that  $a \ne 0$. Recall that, by assumption, we have $b \ne 0$. For any $k=3,\ldots,n$, looking at the coefficients of the  monomials $x_0x_k^2$ and $x_0^2x_k$ in $\overline{F_1}$, we obtain the equations
$$
3am_{0k} + bm_{1k}=0 \ \mbox{ and } \ 3am_{0k}^2+2bm_{0k}m_{1k} =0
$$
which imply that $m_{0k}=m_{1k}=0$ for any $k=3,\dots,n$. 

We may write
$$
R(x_1,x_2)=c_1x_1^3 + c_2x_1^2x_2+c_3x_1x_2^2+c_4x_2^3.
$$
for some $c_1,\dots,c_4\in \mathbb Q$. Looking at the coefficients of the monomials $x_1^2x_k$, $x_1x_k^2$ and $x_2^2x_k$ in $\overline{F_1}$ we see that:
$$
c_2m_{2k}=0 \quad c_3m_{2k}^2=0\quad\text{and} \quad c_4m_{2k}=0.
$$
Since $F$ is a non-degenerate cubic form, it follows that $m_{2k}=0$ for $k=3,\dots,n$. Thus, the claim follows.

\medskip

Assume now that $a=0$. Then, looking at the coefficients of $x_0x^2_k$ and $x_0x_1x_k$, we obtain  $m_{0k}=m_{1k}=0$ for $k=3,\ldots,n$. Thus,  as in the previous case, the claim follows.
\end{proof}

\medskip

\begin{proposition}\label{p_}
Let $F\in \mathbb Z[x_0,\dots,x_n]$ be a non-degenerate cubic form in reduced form:
$$
F(x_0, \ldots, x_n)=ax_0^3 + bx_0^2x_1 + G(x_1,\ldots,x_n)
$$
where $G\in \mathbb Z[x_1,\dots,x_n]_3$. Assume that $\Delta_F \ne 0$.
Let $C \subseteq V_F$ be an irreducible component of positive dimension such that 
$$p=[1,0,\ldots,0] \in C, \quad C \not\subseteq
\{F=0\}\quad \text{and}\quad C \not\subseteq\{x_1=0\}.$$

Then there are finitely many triples 
$$(a_i,b_i,G_i)\in \mathbb Z\times \mathbb Z\times \mathbb Z[x_1,\dots,x_n]_3\qquad i=1,\dots,k$$
such that for all $T\in \SL(n+1,\ZZ)$ such that  $[1,0,\ldots,0] \in T(C)$ 
and $T\cdot F$ is in reduced form, we have that $T\cdot F$ is equivalent to 
$(a_i,(b_i,0),G_i)$ over $\ZZ$ for some $i\in \{1,\dots,k\}$ (cf. Definition \ref{d_reduced}). 
\end{proposition}

\begin{proof}
Suppose not. Then there exist an infinite sequence $T_i\in \SL(n+1,\ZZ)$ with $i=1,2,\dots$ such that  $[1,0,\ldots,0] \in T_i(C)$, $T_i\cdot F$ is in reduced form 
and $T_i\cdot F$ and $T_j\cdot F$ are not equivalent over $\ZZ$ for any $i\neq j$.

 Lemma \ref{l_fp} implies that the set $\{T_i^{-1}([1,0,\dots,0])\}$ is infinite. 
In particular, $C$ admits infinitely many rational points. By Proposition \ref{p_cubic}, we have that $b \ne 0$, as otherwise $p\in W_F$. 

\medskip

We first assume that $C$ is a line. After acting on $(x_1,\ldots,x_n)$ with  $\SL(n,\ZZ)$, we may assume that  $C=\{x_2=x_3=x_4=\ldots=x_n=0\}$ and we may write
$$
F=ax_0^3 + (bx_1+cx_2)x_0^2 + G(x_1,\ldots,x_n)
$$
where $b,c \in \ZZ$, $b \ne 0$ and $G \in \ZZ[x_1,\ldots,x_n]$ is a cubic form. Since reduced forms are considered modulo the action of $\SL(n,\ZZ)$ on $(x_1,\ldots,x_n)$, we may assume that for any $i=1,2,\dots$, the cubic form $F_i=T_i\cdot F$ satisfies the same property, that is 
$$
F_i=a_ix_0^3 + (b_ix_1+c_ix_2)x_0^2 + G_i(x_1,\ldots,x_n)
$$
where $b_i,c_i \in \ZZ$ are such that  $b_i \ne 0$, $G_i \in \ZZ[x_1,\ldots,x_n]_3$ and $T_i(C)=\{x_2=x_3=x_4=\ldots=x_n=0\}$.

Fix $i$ and let $T_i=(t_{hk})_{h,k=0,\ldots,n}$. 
Since $\{x_2=x_3=x_4=\ldots=x_n=0\}$ is fixed by $T_i$ we have $t_{hk}=0$ for $h=2,\ldots,n$ and $k=0,1$. Since $\det T_i=1$, we may assume $\det (t_{h,k})_{h,k=0,1}=1$.

We may find $M, M_i \in \SL(n,\mathbb Q)$  as in  Corollary \ref{c_Cline},  such that
$$
\hat{F}=M \cdot F=ax_0^3+ bx_0^2x_1+dx_1^3 + H(x_2,\ldots,x_n)
$$
and
$$
\hat{F_i}=M_i \cdot F_i = a_ix_0^3+ b_ix_0^2x_1+d_ix_1^3 + H_i(x_2,\ldots,x_n)
$$
where $d,d_i \in \ZZ$ and $H, H_i \in \mathbb Q[x_2,\ldots,x_n]$ are cubic forms.

In addition, if $\hat{T_i}= (\hat t_{hk})_{h,k=0,\dots,n}= M_i \cdot T_i \cdot M^{-1}$, we have that $\hat{T_i}\cdot \hat{F} = \hat F_i$. Let 
$$
U_i:=(\hat t_{hk})_{h,k=0,1}.
$$

Note that, by Corollary \ref{c_Cline}, it follows that $\hat t_{hk}=0$ for $h=2,\ldots,n$ and $k=0,1$ and $U_i \in \SL(2,\ZZ)$. Let
$$
F'=\hat{F}_{|C}=ax_0^3+ bx_0^2x_1+dx_1^3  \mbox{ and } F_i'=\hat F_{i|C}=a_ix_0^3+ b_ix_0^2x_1+d_ix_1^3.
$$
Then  $F',F'_i\in \mathbb Z[x_0,x_1]$ are binary cubic forms such that  $U_i\cdot F'= F_i'$. In particular $\Delta_{F'}=\Delta_{F'_i}\neq 0$
as otherwise  the hypersurface  $\{\hat F=0\} \subseteq\mathbb P^n$ would be  singular and $\Delta_F=\Delta_{\hat F}=0$, which contradicts the assumption on $F$. 
Thus, by Proposition \ref{p_binary} we may assume that 
$$a_i=a \quad b_i=b\quad\text{and} \quad d_i=d\quad\text{for}\quad i=1,2,\dots.$$

On the other hand, by Lemma \ref{l_G}, for each $i=1,2,\dots$
 there 
exists $P_i\in \SL(n-1,\mathbb Q)$ such that $H_i=P_i\cdot H$. Since the hyperplane $\{x_0=0\}$ is invariant with respect to $M_i$, there exist $M,M'_i\in \SL(m,\mathbb Q)$ such that if 
$$
H'(x_1,\dots,x_n)=dx_1^3+H(x_2,\dots,x_n)
$$
and
$$
H'_i(x_1,\dots,x_n)=dx_1^3 + H_i(x_2,\dots,x_n),
$$ 
then $M'\cdot G=H'$ and $M'_i\cdot G=H'_i$ for $i=1,2,\dots$. 
Thus, there exist $P'_i\in \SL(n,\mathbb Q)$ such that $G_i=P'_i\cdot G$ for all $i=1,2,\dots$.
By Jordan's theorem \ref{t_jordan}, it follows that, after possibly taking a subsequence, 
the reduced forms $F_1,F_2,\dots$ are equivalent over $\mathbb Z$. Thus, we obtain  a contradiction.

\medskip

Assume now that $C$ is not a line. Theorem \ref{t_plane} implies that $C$ spans a plane $\Pi$. After acting on $(x_1,\ldots,x_n)$ with  $\SL(n,\ZZ)$, we may assume $\Pi=\{x_3=x_4=\ldots=x_n=0\}$ and we may write
$$
F=ax_0^3 + x_0^2(bx_1+cx_3) + G(x_1,\ldots,x_n)
$$
where $b,c \in \ZZ$, $b \ne 0$ and $G \in \ZZ[x_1,\ldots,x_n]$ is  a cubic form.

Since reduced forms are considered modulo the action of $\SL(n,\ZZ)$ on $(x_1,\ldots,x_n)$, we may assume that this holds for any $i=1,2,\dots$, the cubic form $F_i=T_i\cdot F$ satisfies the same property, that is 
$$
F_i=a_ix_0^3 + x_0^2(b_ix_1+c_ix_3) + G_i(x_1,\ldots,x_n)
$$
where $b_i,c_i\in \mathbb Z$ are such that $b_i\neq 0$, $G_i\in \mathbb Z[x_1,\dots,x_n]_3$ and $T_i(C) \subseteq\Pi=\{ x_3=x_4=\ldots=x_n=0 \}$.

 Fix $i=1,2,\dots$ and let $T_i=(t_{hk})_{h,k=0,\ldots,n}$.   Since $\Pi=\{x_3=\ldots=x_n=0\}$ is fixed by $T_i$ we have $t_{hk}=0$ for $h=3,\ldots,n$ and $k=0,1,2$. Since $\det T_i=1$, we may assume $\det (t_{h,k})_{h,k=0,1,2}=1$.

By Corollary \ref{c_Ccubic}, we may find $M, M_i \in \SL(n,\mathbb Q)$   such that

$$
\hat{F}=M \cdot F=ax_0^3+ bx_0^2x_1+R(x_1,x_2) + H(x_3,\ldots,x_n)
$$
and
$$
\hat F_i=M_i \cdot F_i = a_ix_0^3+ b_ix_0^2x_1+R_i(x_1,x_2) + H_i(x_3,\ldots,x_n)
$$
where  $R,R_i \in \ZZ[x_1,x_2]$ and $H,H_i \in \mathbb Q[x_3,\ldots,x_n]$ are cubic forms. In addition, if $\hat T_i=(\hat t_{hk})_{h,k=0,\dots,n}= M_i \cdot T_i \cdot M^{-1}$, we have that $\hat{T_i}\cdot \hat{F} = \hat F_i$. Let 
$$
U_i:=(\hat t_{hk})_{h,k=0,1,2}.
$$

Note that, by Corollary \ref{c_Ccubic}, it follows that $\hat t_{hk}=0$ for $h=3,\ldots,n$ and $k=0,1,2$ and $U_i \in \SL(3,\ZZ)$. Let
$$
F'=\hat{F}_{|\Pi}=ax_0^3+ bx_0^2x_1+R(x_1,x_2) \ \mbox{ and } {F_i}'=\hat F_{i|\Pi}=a_ix_0^3+ b_ix_0^2x_1+R_i(x_1,x_2).
$$
Then $F',F_i\in \mathbb Z[x_0,x_1,x_2]$ are ternary cubic forms such that $U_i\cdot 
F'= F_i'$. In particular $\Delta_{F'}=\Delta_{F'_i}\neq 0$, as otherwise the hypersurface $\{ \hat F=0\}\subseteq \mathbb P^n$ would be singular and $\Delta_F=\Delta_{\hat F}=0$, which contradicts the assumption on $F$. 
Thus, by  Proposition \ref{p_b2=3} we may assume that  $a_i,b_i$ and $R_i$ do not depend on $i=1,2,\dots$.

As in the previous case, we obtain that, after possibly taking a subsequence,  $F_1,F_2,\dots$ are equivalent over $\mathbb Z$, a contradiction. 
\end{proof}

We are now ready to prove Theorem \ref{t_b2=n}:

\begin{proof}[Proof of Theorem \ref{t_b2=n}]
We may assume that $F$ is in reduced form 
$$
F=ax_0^3 + bx_0^2x_1 + G
$$
where $a, b \in \ZZ$ and $G \in \Z[x_1,\ldots,x_n]_3$.

We assume that there exist $T_i \in \SL(n+1,\ZZ)$, with $i=1,2,\dots$ such that  $F_i=T_i \cdot F$ is in reduced form $(a_i,B_i,G_i)$ for some 
$a_i\in  \mathbb Z$, $B_i\in \mathbb Z^n$ and $G_i\in \mathbb Z[x_1,\dots,x_n]_3$ and  $F_i$ and $F_j$ are not equivalent over $\mathbb Z$ for any $i\neq j$. Acting on $(x_1,\ldots,x_n)$ with $\SL(n,\ZZ)$ we may assume that $B_i=(b_i,0\ldots,0)$, for some $b_i\in \mathbb Z$.  Let $p=[1,0,\dots,0]$ and let $C_1,\dots,C_k \subseteq V_F$ be all the  irreducible components. Then, after possibly replacing $p$ by $T_j(p)$ for some $j$, we may assume that $p,T_i(p)\in C=C_1$ for all $i$ (possibly passing to an infinite subsequence). 
Lemma \ref{l_fp} implies that $C$ is of positive dimension. 

Since by assumption $\Delta_F\neq 0$,   Corollary \ref{c_Csub} implies that  
$$C\nsubseteq \{x_1=0\} \quad\text{and} \quad C\nsubseteq \{F=0\}.$$
Thus, Proposition \ref{p_} implies a contradiction.
\end{proof}


We conclude the section proving a finiteness result on a special class of reduced forms. The result will be used in \S \ref{s_dc}.

\begin{proposition}\label{p_b=0}
Let $F\in \mathbb Z[x_0,\dots,x_n]$ be a non-degenerate cubic form such that $\Delta_F \ne 0$. Fix an integer $r \ne 0$. 
Then there are finitely many pairs 
$$(a_i,G_i)\in \mathbb Z\times \mathbb Z[x_1,\dots,x_n]_3\qquad i=1,\dots,k$$
such that for all $T\in \GL(n+1,\mathbb Z)$ such that $\det T=r$ and $T\cdot F$ is in reduced form, we have that $T\cdot F$ is equivalent to $(a_i,0,G_i)$ over $\mathbb Z$ for some $i \in \{1,\dots,k\}$ (cf. Definition \ref{d_reduced}). Moreover $\Delta_{G_i}\neq 0$ for all $i=1,\dots,k$ 
\end{proposition}

\begin{proof}
Suppose not. Then there exist infinitely many $T_1,T_2,\dots\in \GL(n+1,\mathbb Z)$ such that $\det T_i=r$, $T_i\cdot F=(a_i,0,G_i)$ is in reduced form for each $i$ and $T_i$ and $T_j$ are not equivalent over $\Z$ for each $i\neq j$. 
We denote $S_{i,j}=T_i^{-1}T_j$. Note that 
$T_i([1,0,\dots,0])\in W_F$ for all $i$. Thus,  by Proposition \ref{p_cubic}  we may assume that $[1,0,\dots,0]$ is fixed by $S_{i,j}$ for each $i,j$. 
It follows easily that if $S_{i,j}=(s_{hk})$ then $s_{h0}=s_{0k}=0$ for any $h,k=1,\dots,n$. 

Since $\det T_i=r$, it follows that the denominators of the coefficients of $S_{i,j}$ are bounded and since $\det S_{i,j}=1$, it follows that $s_{0,0}$ is bounded and in particular there exist $i\neq j$ such that $T_i\cdot F$ is equivalent to $T_j\cdot F$ over $\mathbb Z$.

Finally, Lemma \ref{l_disc} implies that, for each $i$ we have $\Delta_{G_i}\neq 0$.
\end{proof}

\medskip

\section{Proof of the main results}

\subsection{Proof of Theorem \ref{t_vol}}\label{s_vol}
Let $X$ be a smooth projective threefold of general type. In this section we prove Theorem \ref{t_vol}, i.e. we show that the volume of $X$ (cf. definition \ref{def_vol}) is bounded by a constant which depends only on the topological  Betti numbers of $X$.

%
%
%

\begin{proof}[Proof of Theorem \ref{t_vol}]
We may assume that $X$ is of general type, as otherwise $\vol(X)=0$. Let $\rmap X.Y.$ be a minimal model of $X$. Then $Y$ has only terminal singularities, and in particular it is smooth outside a finite number of points. In addition, 
$$\vol(X,K_X)=\vol(Y,K_Y)=K_Y^3.$$
Theorem \ref{t_RR} implies that 
$$\chi(Y,\ring Y.)= \frac {1} {24} (-K_Y\cdot c_2(Y)+ e)$$
where 
$$e=\sum_{p_\alpha} \left (r(p_\alpha)-\frac 1 {r(p_\alpha)}\right ),$$
and the sum runs over all the baskets $\mathcal B(Y,p)$ of singularities of $Y$.  
Note that $e\le \Xi(Y)$. 
Thus,
$$\begin{aligned}
\vol(X,K_X)=K_Y^3&\le   3 K_Y\cdot c_2(Y)\\
&= 3 (-24\chi(Y,\ring Y.)+e)\\
&=3(24(-h^{0,0}(X) + h^{1,0}(X) - h^{2,0}(X) + h^{3,0}(X)) + e) \\
&\le 3(12b_3(X)+\Xi(Y)),
\end{aligned}
$$
where the first inequality follows from Theorem \ref{t_BMY} and the second inequality follows from the fact that $h^{1,0} (X)\le h^{2,1}(X)$ by Hard Lefschetz and $h^{2,1}(X)+h^{3,0}(X) \le b_3(X)/2$ by Hodge decomposition.

Thus,  Proposition \ref{p_2b2} implies the claim. 
\end{proof}

Two immediate applications of Theorem \ref{t_vol} are the following corollaries.

\begin{corollary}\label{c_volume}
The volume only takes finitely many values on  the set of three dimensional projective varieties  with a fixed underlying 6-manifold. 
\end{corollary}
\begin{proof} 
Let $X$ be a smooth projective threefold.
The volume $\vol(X,K_X)$ is a rational number whose denominator is bounded by the cube of the index of a minimal model of $X$. By Lemma \ref{l_xi}, the index of any minimal model of $X$ is less than or equal to $4 \cdot \Xi(X)$. The claim follows now from Proposition \ref{p_2b2} and Theorem \ref{t_vol}. 
\end{proof}

\begin{corollary}\label{c_boundedness}
The family of all smooth projective threefolds of general type with bounded Betti numbers is birationally bounded.
\end{corollary}
\begin{proof}
By \cite[Cor. 1.2]{HM06} we know that the family all smooth projective threefolds of general type with bounded volume is birationally bounded. The result follows then from Theorem \ref{t_vol}.
\end{proof}

\subsection{Divisorial contractions}\label{s_dc}

Let $Y$ be a $\Q$-factorial projective threefold and let $f\colon Y \to X$ be an elementary $K_Y$-negative birational contraction. 
By Lemma \ref{l_betti}, we have that $b_2(Y)-b_2(X)=1$. Let $\{\gamma_1, \ldots, \gamma_b\}$ be a basis of  $\bar{H}^2(X,\Z)$ and let $\beta_i=f^*\gamma_i$. 

If $f$ is a divisorial contraction, then we have a natural choice for a class $\alpha \in \bar{H}^2(Y,\mathbb Z)$ such that $\{\alpha,\beta_1,\dots,\beta_b\}$ is a basis of $\bar H^2(Y,\mathbb Q)$ . Indeed,  we can choose $\alpha= c_1(rE)$, where $E$ is the exceptional divisor,  and $r$ is the smallest positive integer such that $rE$ is Cartier. 

 If $f$ is a contraction to a point,  by the projection formula we get
$$
\alpha \cdot \beta_i \cdot \beta_j = 0
$$
and
$$
\alpha^2 \cdot \beta_i=0
$$
for any $i,j=1,\dots,b$. On the other hand, in general , we do not have an isomorphism
$$
\bar{H}^2(X,\Z)=\Z\langle \alpha, \beta_1, \ldots, \beta_b  \rangle
$$
as the following example shows.

\begin{example}
Let $Z=\mathbb P^2$ and consider the $\mathbb P^1$-bundle $$Y=\mathbb P(\ring Z.\oplus \ring Z. (2))$$ over $Z$ with induced morphism $\pi\colon Y\to Z$. Then there exists a birational morphism $f\colon Y\to X$ which contracts a section $E$ of $\pi$ into a point. In particular, $X$ is the cone over $\mathbb P^2$ associated to $\ring Z.(2)$. Note that  $X$ is terminal and $\mathbb Q$-factorial and $K_Y=f^*K_X+1/2 E$.  Let $\ell$ be a line in $Z$ and let $F=\pi^*\ell$. Then $\{E,F\}$ is a basis of $\bar{H}^2(Y,\Z)$. On the other hand, $F'=f_*F$ is not Cartier and therefore it is  not an element of $\bar{H}^2(X,\mathbb Z)$, while $2F'$ is a generator of $\bar{H}^2(X,\mathbb Z)$. 
\end{example}

\bigskip

Given a divisorial contraction to a point  $f\colon Y\to X$ between terminal threefolds, our goal is to first bound the 
difference  $K^3_Y-K^3_X$ and then compute the cubic form $F_X$ associated to $X$ from the cubic form $F_Y$ associated to $Y$. We begin with the following:

\begin{proposition}\label{p_ctp} 
Let $X_0$ be a smooth projective threefold and let 
$$
X_0 \dto X_1 \dto \ldots \dto X_{k-1} \dto X_k
$$ 
be a sequence of steps of the minimal model program for $X_0$. 
Assume that $$f\colon Y=X_{k-1}\to X=X_k$$ is a divisorial contraction to a point $p\in X$.

Then
$$
0<K_{Y}^3-K_{X}^3 \le 2^{10} b_2^{2},
$$
where $b_2=b_2(X_0)$ is the second Betti number of $X_0$.
\end{proposition}
\begin{proof}
Let $E$ be the exceptional divisor of $f$ and let $a=a(E,X)$ be the discrepancy of $f$ along $E$. Since $X$ is terminal, we have that $a>0$. Since  $K_Y^3 - K_X^3= a^3E^3$, it is enough to bound $a^3E^3$. 
The possible values of $aE^3$ are listed in Table 1 and 2 of \cite{Kawhigher}. In particular, we have
$$
0<aE^3 \le 4.
$$	 
Let $\mathcal B(X,p)=\{p_1,\dots,p_k\}$ be the basket of $X$ at $p$, with indices $r_1=r(p_1),\dots,r_k=r(p_k)$ (cf. \S \ref{s_t3}) and let  $R$ be the least common multiple of $r_1,\dots,r_k$. Then,  \cite[Lemma 2.3]{Kawhigher}  implies that $E^3 \ge 1/R$. Thus,
$$
0 < \left(aE\right)^3 \le \frac {64}{(E^3)^2} \le 64 R^2.
$$	

Let $\Xi=\Xi(X,p)\le \Xi(X)$. Then Lemma \ref{l_xi} implies that 
$$
R \le 4\cdot \Xi
$$
and Proposition \ref{p_2b2} implies
$$
\left(aE\right)^3 \le 2^{10} b_2^{2}.
$$
Thus, the claim follows. 
\end{proof}

We now study  the behaviour of the cubic form associated to a terminal threefold, under a divisorial contraction to a point. We begin with the following elementary fact:

\begin{lemma}\label{l_lattice}
Let $A$ be a  maximal rank submodule of $\mathbb Z^m$ and let $r$ be a positive integer.
Assume that for any $b\in \Z^m$ we have that $r\cdot b\in A$.
Let $T\in \mathcal M(m,\Z)$ be a matrix whose columns form a basis of $A$. 

Then $0<|\det T| \le r^m$.
\end{lemma}
\begin{proof} 
By assumption, there exists $X\in \mathcal M (m,\Z)$ such that $T\cdot X=rI_m$, where $I_m\in \SL(m,\mathbb Z)$ is the identity matrix. Thus, $\det T$ divides $r^m$ and the claim follows. 
\end{proof}

\begin{lemma}\label{l_torsion}
Let $X$ and $Y$ be $\mathbb Q$-factorial  projective threefolds with terminal singularities and let $f\colon Y\to X$ be a divisorial contraction onto a point $x \in X$ with exceptional divisor $E$. 

Then $\pi_1(E)=1$ and, in particular, $H^2(E,\Z)$ is torsion-free. 
\end{lemma}

\begin{proof}
Let $U$ be an analytic neighborhood of $x$ such that $U$ retracts to $x$ and consider the morphism $f_U \colon  V=f^{-1}(U) \to U$.  Then, \cite[Theorem 7.8]{Kollar93} implies  that $\pi_1(V)=\pi_1(U)=1$.  
Since $V$ retracts to $E$, it follows that $\pi_1(E)=1$. 

The universal coefficient theorem implies that $H^2(E,\Z)$ is torsion free. 
\end{proof}

Thus, we  have: 
%
%

\begin{proposition}\label{p_cubicpoints}
Let $X$ and $Y$ be $\mathbb Q$-factorial  projective threefolds with terminal singularities and let $f\colon Y\to X$ be a divisorial contraction onto a point with exceptional divisor $E$. 
Let $\alpha\in \bar H^2(Y,\Z)$ be a generator of the ray $\mathbb R_{>0}[E]$ in $N^1(Y)\otimes \mathbb R$.  Let $n=b_2(Y)$ and let $\alpha,\alpha_2,\dots,\alpha_n$ be a basis of $\bar H^2(Y,\Z)$. Let $r=|\alpha^3|$.

Then there exists $T\in \mathcal M(n,\mathbb Z)$ such that $0<|\det T|\le r^n$ and $\alpha, T(\alpha_2), \ldots, T(\alpha_n)$ is a basis of the submodule of $\bar H^2(Y,\Z)$ spanned by $f^* \bar H^2(X,\Z)$ and $\alpha$.  

In particular, it follows that 
$$
T\cdot F_Y= ax_0^3  +  F_X(x_1,\ldots,x_n),
$$
where $a=\alpha^3$.
\end{proposition}
\begin{proof}
Fix an isomorphism $\bar H^2(Y,\Z)\simeq \mathbb Z^n$ and consider the submodule $A$ of $\Z^n$ spanned by $f^* \bar H^2(X,\Z)$ and $\alpha$.  Let $\beta\in \bar H^2(Y,\Z)$. Then there exist integers $c,b$ with $|b| \le r$ such that
$$
(c \alpha + b \beta) . \alpha^2=0.
$$

Set $\gamma=c \alpha + b \beta$. As in the proof of Lemma \ref{l_betti}, it follows that $R^1f_*\mathbb Z=0$ and, in particular, $H^1(E,\Z)=0$. Thus,  as in Lemma \ref{l_cohomology}, we get the exact sequence

$$
0 \to f^*\bar H^2(X,\Z) \to \bar H^2(Y,\Z) \xrightarrow{p} H^2(E,\Z). 
$$

Possibly passing to a desingularization, we can apply \cite[Proposition 12.1.6]{KM92} to obtain  that $p(E)$ is a multiple of  $p(\gamma)$  in $H^2(E,\Q)$. Since $\gamma.\alpha^2=0$, it follows that $p(\gamma)$ is a torsion element of $H^2(E,\Z)$, which implies that $p(\gamma)=0$ by Lemma \ref{l_torsion} and so  $\gamma \in f^*\bar H^2(X,\Z)$.

 Thus, $b\beta\in A$ and Lemma \ref{l_lattice} implies the claim. 
\end{proof}

%
%
%
%
%
%
%
%
%

\medskip

We now consider divisorial contraction to a smooth curve. 
We begin with the following well known result (e.g. see \cite[Prop. 14]{OV95}):

\begin{proposition}\label{l_F_Y}
Let $X$ be a $\mathbb Q$-factorial projective threefold and let $C$ be a smooth curve of genus $g$  contained in the smooth locus of $X$. Let
$f\colon Y \to	X$ be the blow-up of $X$ along $C$ and let  $\alpha=c_1(E)$.  

Then $H^2(Y,\Z) \cong \Z[\alpha] \bigoplus H^2(X,\Z)$ and
$$
K_Y^3-K_X^3= -2K_X.C+2-2g=-2E^3 +6-6g.
$$
In particular, if $\beta_1,\dots,\beta_n$ is a basis of $H^2(X,\mathbb Z)$, then 
$\alpha,f^*\beta_1,\dots,f^*\beta_n$ is a basis of $H^2(Y,\Z)$ and with respect to these basis we have:
$$
F_Y(x_0,\ldots,x_n)= ax_0^3  + 3x_0^2(\sum_{i=1}^n b_i x_i) + F_X(x_1,\ldots,x_n),
$$
where $a=\alpha^3$ and $b_i= -\beta_i.C$.
\end{proposition}

\subsection{Proof of Theorem \ref{t_global}} We can finally prove our main result.

\begin{proof}[Proof of Theorem \ref{t_global}]
 
If $f$ is a divisorial contraction to a point, then $(1)$ is the content of Proposition \ref{p_ctp}. Assume hance that $f$ contracts a divisor $E$ to a smooth curve $C$. Then $E$ is a $\mathbb P^1$-bundle over $C$ and, in particular, if $g$ is the genus of $C$ then $b_3(E)=2g$. Thus, 
by Lemma \ref{l_cohomology} and Lemma \ref{l_betti} and since $E$ and $C$ are smooth, we have that
$$
b_3(Y)-b_3(X)=Ib_3(Y)-Ib_3(X)=2g.
$$

Moreover, considering the cubic form $F_Y$ associated to $Y$ and applying Proposition \ref{l_F_Y},  we have that $|E^3| \le S_Y$.  Hence
\begin{align*}
|K_Y^3-K_X^3| &=|-2E^3 +6-6g| \\ &\le 2 S_Y+6(b_3(Y)+1).
\end{align*}
This finishes the proof of  $(1)$.

\medskip
We now prove $(2)$.

Let us first assume that $f$ is a divisorial contraction to a point with exceptional divisor $E$. Let $\alpha \in H^2(Y,\Z)$ be a generator of the ray $\mathbb R_{>0}[E]$. By Propostion \ref{p_cubicpoints}, $\alpha$ is a point of rank 1 for the Hessian of the cubic form $F_Y$. Then, by Proposition \ref{p_cubic}, $\alpha$ is determined up to finite ambiguity by $F_Y$ and so it is $r=\alpha^3$. By Proposition \ref{p_b=0}, there are finitely many pairs
$$(a_i,G_i)\in \mathbb Z\times \mathbb Z[x_1,\dots,x_n]_3\qquad i=1,\dots,k$$
such that for all $T\in \mathcal {M}(n+1,\mathbb Z)$ such that $0 <|\det T| \le r^n$ and $T\cdot F$ is in reduced form, we have that $T\cdot F$ is equivalent to $(a_i,0,G_i)$ over $\mathbb Z$ for some $i \in \{1,\dots,k\}$. By Proposition \ref{p_cubicpoints},  there exists $T\in \mathcal{M}(n+1,\mathbb Z)$ such that $0<|\det T|\le r^n$ and 
$T\cdot F$ is in reduced form $(a,0,F_Y)$, where $a=\alpha^3$. Thus, there exists $M\in  \SL(n,\mathbb Z)$ such that 
$F_Y=M\cdot G_i$ for some $i\in \{1,\dots,k\}$.

Let us assume now that $f$ is a divisorial contraction to a smooth curve. By Theorem  \ref{t_b2=n}, there exist finitely many tripes 
$$(a_i,B_i,G_i)\in \mathbb Z\times \mathbb Z^{n}\times \mathbb Z[x_0,\dots,x_{n}]_3\qquad i=1,\dots,k$$
such that any reduced triple associated to $F$ is equivalent to $(a_i,B_i,G_i)$ over $\mathbb Z$ for some $i\in\{1,\dots,k\}$. By Proposition \ref{l_F_Y}, 
there exist $a\in \mathbb Z$ and $B\in \mathbb Z^{n}$ such that $(a,B,F_Y)$ is a reduced triple associated to $F$. Thus, there exists $M\in  \SL(n,\mathbb Z)$ such that 
$F_Y=M\cdot G_i$ for some $i\in \{1,\dots,k\}$

\end{proof}

\begin{proof}[Proof of Corollary \ref{c_top}]
This is a simple iteration of Point $(2)$ of Theorem \ref{t_global}, keeping in mind that if $g: W \to Z$ is a step of an MMP as in the statement and $\Delta_{F_W} \ne 0$, then also $\Delta_{F_Z} \ne 0$ (this follows combining together Proposition \ref{l_F_Y} and Proposition \ref{p_cubicpoints} with Lemma \ref{l_disc}). 
\end{proof}

\begin{proof}[Proof of Corollary \ref{c_main}]
Let 
$$
X=X_0 \to X_1 \to \ldots \to X_k
$$
be an MMP for $X$ such that each $f_i\colon X_i \to X_{i+1}$ is a divisorial contraction to a point or to a smooth curve contained  in the smooth locus of $X_{i+1}$.

Denote by $F_{i}$ the cubic form associated to $X_{i}$ and let $S_i=S_{X_i}$ (cf. Definition \ref{d_skansen}). Theorem \ref{t_b2=n} implies that $S_{X_0}<+\infty$.

We proceed by induction on $i=0,\ldots,k$.
Proceeding as in the proof of Theorem \ref{t_global}, by combining together Proposition \ref{l_F_Y}, Proposition \ref{p_cubicpoints}, Proposition \ref{p_b=0} and Theorem \ref{t_b2=n}, it 
follows that, for any $i=0,\ldots,k$,   
$$\Delta_{F_i} \ne 0\qquad\text{and}\qquad S_i < + \infty.$$ 
Moreover, each $S_i$ depends only on $F_X$ and, therefore,  only on the topology of the manifold underlying $X$. 

We define
$$
D_k=6b_2(X)+36b_3(X)
$$ 
and for any $i<k$,  let 
$$
D_i=D_{i+1} + \max\{2^{10} b_2(X)^{2}, 2S_i + 6(Ib_3(X_i)+1)  \}.
$$

We claim that
$$
|K_{X_i}^3| \le D_i
$$
for any $i=0,\ldots,k$.

The proof is by descending induction on $i=k, \ldots, 0$. If $i=k$ the result is exactly Theorem \ref{t_vol}. Assume now that $i < k$ and $|K_{X_{i+1}}^3| \le D_{i+1}$. Then the claim follows by combining Proposition \ref{p_ctp} and Theorem \ref{t_global}.
In particular, we have that $|K_X^3|\le D_0$.  

Finally, we need to show that for any $i=1,\dots,k$, we have that $Ib_3(X_i)\le Ib_3(X_{i-1})$. If $f_{i-1}\colon X_{i-1}\to X_{i}$ is a divisorial contraction to a point, then the claim follows immediately from Lemma \ref{l_cohomology}. On the other hand, if $f_{i-1}\colon X_{i-1}\to X_{i}$ is a divisorial contraction to a smooth curve $C_i\subseteq X_{i}$  with exceptional divisor $E_i$, then $E_i$ is a $\mathbb P^1$-bundle over $C_i$ and if $g(C_i)$ is the genus of $C_i$ then  Lemma \ref{l_cohomology} implies
$$Ib_3(X_{i-1})-Ib_3(X_{i})=Ib_3(E_i)=b_3(E_i)=2g(C_i)\ge 0,$$ 
as claimed.

Thus,  $Ib_3(X_i) \le Ib_3(X)=b_3(X)$ for any $i=1,\ldots,k$ and  the Theorem follows. 
\end{proof}

\bibliographystyle{amsalpha}
\bibliography{Library}
\end{document}